\documentclass{amsart}
\usepackage{amsmath,amssymb,enumerate}

\newtheorem{prop}{Proposition}[section]
\newtheorem{coro}[prop]{Corollary}
\newtheorem{lem}[prop]{Lemma}
\newtheorem{rem}[prop]{Remark}
\newtheorem{exe}[prop]{Example}
\newtheorem{defi}[prop]{Definition}
\newtheorem{theo}[prop]{Theorem}

\newtheorem{theoi}{Theorem}

\newcommand{\R}{\mathbb R}

\newcommand{\C}{\mathbb C}
\newcommand{\N}{\mathbb N}
\newcommand{\V}{\mathbb V}
\newcommand{\Z}{\mathbb Z}
\newcommand{\Or}{\hat{\mathcal O}_{\rho}}
\newcommand{\Os}{\hat{\mathcal O}_{\sigma}}
\newcommand{\Ar}{\underline{\mathrm{Art}}}
\newcommand{\Orn}{\mathcal O_{\rho | n}}
\newcommand{\Otn}{(\hat{\mathcal O}_{T} / \mathfrak{m}^n)}
\newcommand{\Ot}{\hat{\mathcal O}_{T} }
\newcommand{\rep}{R(\Gamma, GL_N)}
\newcommand{\repx}{R(\Gamma, GL(\mathbb{V}_{\rho,x}))}
\newcommand{\repg}{R(\Gamma, G)}
\newcommand{\m}{\mathfrak m}
\newcommand{\jo}{\mathfrak j}
\newcommand{\hdot}{\bullet}

\newcommand{\gu}{{\mathfrak g}}
\newcommand{\hu}{{\mathfrak h}}

\newenvironment{prf}{Proof:}{$\Box$}

\begin{document}

\title[VMHS and deformation theory of a C-VHS]{Variations of Mixed Hodge Structure attached to the deformation theory of a Complex Variation of Hodge Structures}
\author[P. Eyssidieux]{Philippe Eyssidieux}
\address{Institut Fourier, UMR 5582 CNRS-Universit\'e Joseph Fourier
\\ Universit\'e de Grenoble I\\
100 rue des Maths, BP 74, 38402 Saint Martin d'H\`eres Cedex, France}
\email{philippe.eyssidieux@ujf-grenoble.fr}
\urladdr{http://www-fourier.ujf-grenoble.fr/$\sim$eyssi/}
\author[C. Simpson]{Carlos Simpson}
\address{CNRS, Laboratoire J. A. Dieudonn\'e, UMR 6621
\\ Universit\'e de Nice-Sophia Antipolis\\
06108 Nice, Cedex 2, France}
\email{carlos@unice.fr}
\urladdr{http://math.unice.fr/$\sim$carlos/}
\thanks{This research is partially supported
by ANR grants BLAN08-1-309225 (SEDIGA) and BLAN08-3-352054 (G-FIB)} 
\date{February 11, 2009}

\maketitle

Let $X$ be a compact connected K\"ahler manifold, $x\in X$ and $\Gamma=\pi_1(X,x)$. 
Let $\rho: \Gamma \to GL_N(\C)$ be a finite dimensional semisimple representation.
We assume $\rho$ to be the monodromy of a given polarized $\C$-VHS 
$(\mathbb{V}_{\rho}, \mathcal{F}^{\bullet},\overline{ \mathcal {G}}^{\bullet},
  S)$  whose weight is zero. If $\rho$ is not irreducible
then several distinct polarizations could be chosen, we fix one once for all. 
In the introduction, we fix an isomorphism $\mathbb{V}_{\rho,x}\to \C^n$.  

 Then, the Zariski closure of its monodromy group is a reductive subgroup $G\subset GL_N$. Let $\rep$ be the variety of its representations in $GL_N$ \cite{LuMa}. $\rep$ may be viewed as an affine scheme  over $\Z$ but we will only consider it as an affine scheme
over $\C$. The group $GL_N$ acts algebraically on $\rep$ by conjugation and we denote by $\Omega_{\rho}$
the orbit of $\rho$. It is a closed smooth algebraic subvariety and we will consider it as a subscheme of $\rep$ endowing it with its reduced induced structure. 

Denote by $\rep_{\rho}$ the formal local scheme which is the germ at $[\rho]$ of $\rep$. Similarly, denote by $\hat\Omega_{\rho}$ the germ of $\Omega_{\rho}$ at $[\rho]$. $\hat \Omega_{\rho}$ is a closed formal subscheme of $\rep_{\rho}$. 

Let $(\Or,\m)$ be the complete local ring of $[\rho] \in\rep(\C)$ so that $$\rep_{\rho}=Spf(\Or).$$
Let  $\jo \subset \Or$ the ideal defining $\hat\Omega_{\rho}$. $\jo$ is a prime ideal. 
Let $n\in \N^*$ be a positive integer and $\Orn:= \Or/ \m ^n$ be the Artin algebra of the $n-1$-th infinitesimal neighborhood of $[\rho]$.
 For $n\ge k\ge 0$, let $W_{-k} \Orn\simeq \jo^k/\jo^k\cap \m^n$ be the image of $\jo^k$ in $\Orn$.

 Let $T=Spf(\Ot)$ be the formal germ at the origin of the 
homogenous quadratic cone defined by the zeroes of the obstruction map $$\mathrm{obs_2}: S^2H^1(X, \mathrm{End}(\V_{\rho}))\to H^2(X,\mathrm{End}(\V_{\rho})).$$

One can  slightly revisit the statements of \cite{GM} and give the construction of an isomorphism of formal germs:
$$ GM^c: Spf(\Or) \to \hat \Omega_{\rho} \times T
$$
which we call the preferred Goldman-Millson isomorphism (definition \ref{can}). 
The main theorem of \cite{GM} asserts that an abstract isomorphism between these two formal germs exists.
 The original construction is not quite canonical, and we have tried to make the choices made there as explicit as possible. This is the rationale for the length of section 2. It turns out that there is no good reason to give any kind of privilege to this preferred isomorphism. Two other very nice Goldman-Millson isomorphisms 
$GM'$ and $GM''$ can be constructed and the rest of the article can in retrospect be understood as the study their interplay.

 A direct consequence of \cite{GM}
is the construction of a split MHS on $\Orn$ underlying this weight filtration. 
The MHS on $\Or$ obtained by passing to the limit is rather unpleasant since it has infinite dimensional weight subquotients. The choices made in the construction 
and the fact that the MHS is split make it clear that it is not the right one, which was conjectured to exist in \cite{MTS},  but only its $Gr_{W}$. Hence this MHS is only an approximation of the true object 
\footnote{
It should be noted that  this  construction is essentially included in \cite[p. 5]{Pri1}, see also \cite{Pri3}.}.

Since the outcome of this construction does not depend on $x\in X$, something more subtle must 
be done to confirm the conjecture in \cite{MTS} that there should be a (non necessarily split) MHS on $\Orn$ depending explicitely on $x$.  This is done in section 3. 

The complete local ring $\Ot$ is somewhat better
behaved than $\Or$ in certain respects. In particular, the graded Artin local ring $\Otn$  has a canonical split  MHS with finite dimensionnal weight subquotients and weight filtrations given by the powers of the maximal ideal. On the other hand, its mixed Hodge Theory is a bit more messy than that of $\Or$.

Consider the tautological representation $\rho_n: \Gamma \to GL_N\Otn$ induced by the natural 
map $T\to Spf(\Or), \ t\mapsto (GM^c)^{-1}([\rho], t)$ and $\V_n$ the locally free rank $N$ $\Otn$-local system
on $X$ attached to $\rho_n$. Since $\rho_n(\Gamma)$ preserves the decreasing filtration defined at level $k$ by $\m^k \Otn^N \subset \Otn^N$, this filtration induces on $\V_n$  an  increasing filtration indexed by the nonpositive integers
$W_{-n} \subset \ldots \subset W_{-k}\subset W_{-k+1} \subset \ldots \subset W_{0}= \V_n$.

Our main result is the following basic fact:

\begin{theoi}\label{gln}
The filtration $W_{\bullet}$ is the weight filtration of a $\C$-Variation of Mixed Hodge Structures on $X$. 
\end{theoi}

The proof is given in section 3. These VMHS are closely related to -but somewhat different from-
those constructed in \cite{Ha2}, Theorem 13.10 p. 82 \footnote{We have not fully understood this link yet.}. The explicit nature of their 
construction  makes them easier to use in applications. 

In fact, these results admit an easy  generalization developped in section 4. 

\begin{theoi} \label{gengp}
Let $G$ be a reductive algebraic group defined over $\C$. Let $\sigma: \Gamma\to G(\C)$ be a semisimple representation 
whose associated Higgs bundle is a fixed point of the $\C^*$-action on $M_{Dol}(X,G)$ \cite{Sim3}.
 
 Let $T=Spf(\Ot)$ be the formal germ at the origin of the 
homogenous quadratic cone attached to the obstruction map $$S^2H^1(X, \mathrm{ad}_{\sigma})\to H^2(X,\mathrm{ad}_{\sigma}).$$

Let $\repg/\C$ be the affine scheme parametrizing the representations of $\Gamma$ with values in $G$ endowed with 
the action of $G$ by conjugation. Let $\hat\Omega_{\sigma}$ be the formal germ at $[\sigma]$ of the orbit of $\sigma$. Then, there is a preferred isomorphism $GM^c: Spf(\Os) \to \hat\Omega_{\sigma}\times T$. 

For every $n$, $\Os/\m^n$ carries a canonical functorial $\C$-MHS. 

Let $\alpha$ be a rational  representation of $G$ with values in $GL_N$
 and let $\sigma_n:\Gamma\to G\Otn$ be the tautological representation. Denote by $\V_{\alpha,\sigma}$ the  local system in $\Otn$ free modules on $X$
 attached to the representation $\alpha\circ \sigma:\Gamma\to GL_N\Otn$.
  
  The $\C$-local system 
 underlying $\V_{\alpha,\sigma}$
 is the holonomy of a graded polarizable VMHS whose weight filtration is given by $$W_{-k} \V_{\alpha,\sigma}= \m^k.\V_{\alpha,\sigma} \quad k=0,..., n.$$ 
\end{theoi}

\section{Basic definitions in Mixed Hodge Theory}

The following definitions are well-known to the experts. We nevertheless recall them for the reader's convenience. 

\subsection{$\C$-MHS}

\begin{defi}
Let $w\in\Z$ be an integer. 
A complex (finite dimensionnal) Hodge structure (in short, HS) of weight $w$ is a triple 
 $(V, F^{\hdot}, \bar G^{\hdot})$ where $F^{\hdot}$, $\bar G^{\hdot}$ are $w$-opposed 
decreasing biregular filtrations on the (finite dimensionnal) $\C$-vector space $V$,  that is 
 $Gr^p_{F}Gr^{q}_{\bar G} V= 0$ if $p+q\not = w$. 
 
 If $V$ is finite dimensionnal, a hermitian form $S$ on  $V$ polarizes 
 the Hodge Structure if,  using the usual definition $H^{p,q}= F^p\cap \bar G^q$ for $p+q=w$, the decomposition $V=\oplus _{p+q=w} H^{p,q}$ is $S$-orthonormal and
 $(-1)^{p+w} S |_{H^{p,q}}>0$.
 
 A complex (finite dimensionnal) Mixed Hodge Structure (in short, MHS) is a quadruple  $(V, F^{\hdot}, \bar G^{\hdot}, W_{\hdot})$
 where $F^{\hdot}$, $\bar G^{\hdot}$ 
decreasing biregular filtrations and $W_{\hdot}$ is an increasing filtration on the (finite dimensionnal) $\C$-vector space $V$
such that the filtrations induced by $F^{\hdot}, \bar G^{\hdot}$ on $Gr_{k}^W V$ give rise to a complex Hodge Structure. 

A complex (finite dimensionnal) HS defined over the real is a $\R$-vector space $V_{\R}$ with a $\C$-Hodge structure $(V_{\C}, F^{\hdot}, \bar F^{\hdot})$
on its complexification, such that $\bar F^{\hdot}$ is the complex conjugate of $F^{\hdot}$.

A complex (finite dimensionnal) MHS defined over the real is a $\R$-vector space $V_{\R}$ with a $\C$-Mixed Hodge structure $(V_{\C}, W_{\hdot},  F^{\hdot}, \bar F^{\hdot})$
on its complexification, such that $W_{\hdot}$ is the complexification of an increasing 
 filtration of $V_{\R}$ and $\bar F^{\hdot}$ is the complex conjugate of $F^{\hdot}$. 
\end{defi}

\begin{prop} (\cite{H2}) The category of complex Mixed Hodge structures is abelian.
For the usual tensor product and duality functor the category of complex Mixed Hodge structures, it is tannakian. 
$\C$-Hodge structures form a full abelian and tannakian subcategory.
\end{prop}

The following definition is thus natural:
\begin{defi}
A $\C$- Hodge (positively) graded Lie algebra is a complex graded vector space $\oplus_{k\ge 0} L^k$, each $L^k$ being endowed with 
a $\C$-Hodge structure of weight $k$, with a graded Lie algebra bracket respecting the Hodge structure. 

A Hodge Lie Algebra is a real finite dimensional Lie Algebra $\mathfrak g$ such that $\mathfrak g_{\C}$ carries
a $\C$-Hodge structure of weight zero defined over $\R$ respected by the Lie bracket. 
\end{defi}

\begin{exe}
The Lie algebra of a group of Hodge type is a Hodge Lie algebra.  
\end{exe}

\begin{lem}\label{twist}
Let $\mathbb{M}=(V, F^{\hdot}, \bar G^{\hdot}, W_{\hdot})$ be a Mixed Hodge Structure. Let $u\in GL(V)$ such that $u-Id_V \in W^{-1} End(V)$. Then $\mathbb{M}_u=((V, F^{\hdot}, u(\bar G^{\hdot}), W_{\hdot})$ is a MHS such that
$Gr_{\hdot}^W \mathbb{M}_u= Gr_{\hdot}^W \mathbb{M}$ as split Mixed Hodge Structures.
\end{lem}
\begin{prf}
This is immediate since the filtration induced by $ \bar G^{\hdot}$ and $u( \bar G^{\hdot})$ on $Gr_{\hdot}^W V$ co\"incide. Indeed $u$ preserves $W^{\hdot}$ and induces identity on $Gr_{\hdot}^W V$. 
\end{prf}

\begin{lem}\label{4f}
Let $(V, F^{\hdot}, \bar G^{\hdot}, W_{\hdot})$ be a finite dimensional complex vector space with three biregular filtrations having a fourth biregular filtration
$U^{\hdot}$ such that induced filtrations induced by $ F^{\hdot}, \bar G^{\hdot}, W_{\hdot}$ on every $Gr_{U}^{r} V$
give rise to a Mixed Hodge Structure.
Then  $(V, F^{\hdot}, \bar G^{\hdot}, W_{\hdot})$ is itself a Mixed Hodge Structure. 
\end{lem}
\begin{prf}
There should be a reference for the claim in the MHS literature, but here is a short proof. 
It suffices
by induction to prove this for a $2$-step filtration, which comes down to the following statement: if
$U$ is a vector space with three filtrations $W,F,G$, and if $U'\subset U$ is a subspace with quotient denoted
$U'':= U/U'$, if we suppose that the induced filtrations on $U'$ and $U''$ are HMS, then the filtrations on $U$ formed
a MHS. 

In turn, this can be seen by Penacchio's interpretation \cite{Pen}: given three filtrations we get a bundle on ${\mathbb P}^2$,
and they form a MHS if and only if the bundle is $\mu$-semistable of slope $0$.  The subspace $U'$ with its filtrations corresponds to a
locally free subsheaf, and $U''$ is the reflexive sheaf associated to the quotient sheaf. We get a short exact sequence of sheaves outside of codimension $2$,
and in this case if the kernel and cokernel are $\mu$-semistable of slope $0$ then so is the middle bundle. 
\end{prf}

Suppose $A$ is an artinian local ring, with increasing filtration $W$ and decreasing filtrations $F$ and $G$,
all compatible with the algebra structure. Let ${\bf m}$ denote the maximal ideal and let $M^k:= {\bf m}^k$
be the decreasing filtration of $A$ by powers of $n$. Let $n$ be the smallest integer with $M^{n}=0$.

Denote by $V:=M^1/M^2 = Gr^1_M(A)$, which is the dual of the Zariski tangent space of $A$. Note that the filtrations
$W$, $F$ and $G$ induce filtrations given by the same letter on the associated graded pieces $Gr^k_M(A)$, in particular on $V$.
The associated-graded algebra $Gr^{\cdot}_M(A)= \bigoplus _{k=0}^{n -1}Gr^k_M(A)$ is generated by the
piece in degree $1$ which is $V$. This means that we have surjections 
$$
Sym^k(V)\stackrel{\mu ^k}{\rightarrow} Gr^k_M(A) \rightarrow 0.
$$
The source $Sym^k(V)$ has three filtrations obtained from the symmetric product operation applied to the filtrations 
of $V$ whereas the target $Gr^k_M(A)$ has induced filtrations as stated above. The map $\mu ^k$ comes from the algebra structure so it 
preserves the three filtrations.

\begin{prop} \label{mhalg}
Suppose that the above data satisfy the following hypotheses:
\newline
(1)---$Spec(A)$ is the $n$-th neighborhood of the  origin in a quadratic cone;
\newline
(2)---the filtrations $W,F,G$ induce a complex MHS on $V$, which we use also to give a CMHS on $Sym^k(V)$;
\newline
(3)---the kernel $K$ of the map $\mu ^2:Sym^2(V)\rightarrow Gr^2_M(A)$ is a sub-CMHS of $Sym^2(V)$; and 
\newline
(4)---for each $k$, the filtrations induced by $W$, $F$ and $G$ on $Gr^k_M(A)$ are the same
as the quotient filtrations induced by the map $\mu ^k$ from the filtrations on $Sym^k(V)$,
in other words $\mu ^n$ strictly preserves the filtrations.

Then $W,F,G$ induce a CMHS on $A$ and $M^{\cdot}$ is a filtration by sub-CMHS's. 
\end{prop}
\begin{prf}
Using condition (1) we get exact sequences 
$$
K\otimes Sym^{k-2}(V)\rightarrow Sym^k(V)\rightarrow Gr^k_M(A)\rightarrow 0.
$$
By condition (3), $K$ is a sub-CMHS of $Sym^2(V)$, so the map 
$K\otimes Sym^{k-2}(V)\rightarrow Sym^k(V)$ is a morphism of CMHS. Hence, its cokernel is a CMHS,
then condition (4) says that $Gr^k_M(A)$ with its triple of filtrations is equal to this cokernel,
so $Gr^k_M(A)$ is a CMHS. Now, apply lemma \ref{4f}. 

\end{prf}

\subsection{$\C$-VMHS}

\begin{defi}
A $\C$-VHS (polarized complex variation of Hodge structures)
  on $X$ of weight $w\in \Z$ is a 5-tuple $(X,
  \mathbb{V},\mathcal{F}^{\bullet},\overline{ \mathcal {G}}^{\bullet},
  S)$ where:
\begin{enumerate}
\item $\mathbb{V}$ is a local system of finite dimensional
  $\C$-vector spaces, 
\item   $S$ a non degenerate flat sesquilinear pairing on
  $\mathbb{V}$,  
\item $\mathcal{F}^{\bullet} =(\mathcal{F}^p)_{p\in \Z}$ a biregular
  decreasing filtration of $\mathbb{V}\otimes_{\mathbb C}
  \mathcal{O}_X$ 
by locally free coherent analytic sheaves    such that $d' {\mathcal F}^p \subset \mathcal F ^{p-1} \otimes
  \Omega^1_X$, 
\item $ \overline{ \mathcal {G}}^{\bullet} =(\overline{ \mathcal
  {G}}^{q})_{q\in\Z}$ a biregular decreasing filtration of
  $\mathbb{V}\otimes_{\mathbb C} \mathcal{O}_{\bar X}$ 
by locally free coherent antianalytic sheaves such that
 $d''\overline{\mathcal{G}}^p \subset \overline{\mathcal{G}
}^{p-1} \otimes \Omega^1_ {\bar X}$,
\item for every point $x\in X$ the fiber at $x$ $
  (\mathbb{V}_x,\mathcal{F}_x^{\bullet},\overline{ \mathcal {G}}_x^{\bullet})$ is a $\C$-HS
  polarized by 
  $S_x$.
\end{enumerate}
\end{defi}

 This definition is easily seen to be equivalent to that given by
\cite{Sim1}.

\begin{exe} \label{hdga}
Let $\rho:\pi_1(X,x)\to GL(\V_{\rho,x})$ be the monodromy representation underlying a $\C$-VHS and $\V_{\rho}$ be the correpsonding local system, then $ad_{\rho}=\mathrm{End}(\V_{\rho})$ is a local system in  Lie algebras that underlies a $\C$-VHS.

 Then $H^{\hdot}(X, \mathrm{End}(\V_{\rho}))$ is a $\C$-Hodge graded Lie algebra for the bracket obtained by composing the usual cup product
 $$H^{\hdot}(X, \mathrm{End}(\V_{\rho})) \otimes H^{\hdot}(X, \mathrm{End}(\V_{\rho})) \to H^{\hdot}(X, \mathrm{End}(\V_{\rho})\otimes \mathrm{End}(\V_{\rho}))$$ with the
 cohomology operation $$H^{\hdot}(X, \mathrm{End}(\V_{\rho})\otimes \mathrm{End}(\V_{\rho}))\to H^{\hdot}(X,  \mathrm{End}(\V_{\rho}))$$ induced by the Lie bracket 
 $ \mathrm{End}(\V_{\rho})\otimes \mathrm{End}(\V_{\rho})\to \mathrm{End}(\V_{\rho})$. 
\end{exe}

The $\C$-Hodge structure on the cohomology of a $\C$-VHS can be constructed by a straightforward adaptation of the Deligne-Zucker argument \cite{Zuc},  which is written for the real case. 

The following definition is a slight generalisation of the definition in \cite{U}. 

\begin{defi}
A $\C$-VMHS
on $X$ is a 6-tuple $(X,
\mathbb{V},\mathbb{W}_{\bullet},\mathcal{F}^{\bullet},\overline{
\mathcal {G}}^{\bullet}, (S_k)_{k\in\Z})$ where:
\begin{enumerate}
\item $\mathbb{V}$ is a local system of finite dimensional
  $\C$-vector spaces, 
\item $\mathbb{W}_{\bullet}= (\mathbb{W}_k) _{k\in \Z}$ is a
  decreasing filtration of $\mathbb{V}$ by local subsystems,   
\item $\mathcal{F}^{\bullet} =(\mathcal{F}^p)_{p\in \Z}$ a biregular
  decreasing filtration of $\mathbb{V}\otimes_{\mathbb C}
  \mathcal{O}_X$ 
by locally free coherent analytic sheaves  such that $d' {\mathcal F}^p \subset \mathcal F ^{p-1} \otimes
  \Omega^1_X$, 
\item $ \overline{ \mathcal {G}}^{\bullet} =(\overline{ \mathcal
  {G}}^{q})_{q\in\Z}$ a biregular decreasing filtration of
  $\mathbb{V}\otimes_{\mathbb C} \mathcal{O}_{\bar X}$ 
by locally free coherent antianalytic sheaves such that $d''\overline{\mathcal{G}}^p \subset \overline{\mathcal{G}
}^{p-1} \otimes \Omega^1_ {\bar X}$, 
\item  $\forall x\in X$ the stalk $(\V_x,  \mathbb{W}_{\bullet,x},\mathcal{F}^{\bullet}_x,\overline{
\mathcal {G}}^{\bullet}_x)$ is a $\C$-MHS,
\item  $S_k$ is flat sesquilinear non degenerate pairing on
  $Gr^{\mathbb{W}}_k \mathbb{V}$, 
\item $(X, Gr^{\mathbb{W}}_k {\mathbb V},   \mathcal{F}^{\bullet} \cap
  Gr^{\mathbb{W}}_k{\mathbb V}\otimes_{\C} O_X,
  \overline{\mathcal{G}}^{\bullet}\cap Gr^{\mathbb{W}}_k{\mathbb
    V}\otimes_{\C} O_{\bar X}, S_k)$ is a $\C$-VHS.  
\end{enumerate}
\end{defi}

A $\C$-VMHS is uniquely determined by its monodromy (as a
$\mathbb{W}$-filtered representation of $\pi_1(X,x)$) and the MHS
$(\mathbb{V}_x,\mathbb{W}_x^{\bullet},\mathcal{F}_x^{\bullet},\overline{
\mathcal {G}}_x^{\bullet} )$ (rigidity theorem,
cf. \cite[p.85,(1.7)c]{HZ}) and the references therein) .

To be consistant with earlier terminology, we could also have called  $\C$-VHS polarized $\C$-VHS (resp. $\C$-VMHS  graded polarized $\C$-VMHS ).

\begin{lem} 
 Let $E^{\hdot}=E^{\hdot}(X, \mathrm{End}(\V_{\rho}))$ be the $C^{\infty}$-De Rham complex of the $\C$-VHS 
attached to $ad_ {\rho}=\mathrm{End}(\V_{\rho})$. It is endowed with the usual Hodge filtrations and the usual Lie bracket preserves the Hodge filtrations. We have the familiar conditions from \cite{H3}:

\begin{itemize}
 \item The differential of $E^{\hdot}$ is strictly compatible with 
the two filtrations induced by $F^{\hdot}$ and $\bar G^{\hdot}$.
\item The induced filtrations on $H^k(X,\mathrm{End}(\V_{\rho}))$ give a $\C$-Hodge Structure of weight $k$. 
\end{itemize}

The Lie bracket induces on the cohomology $H^{\hdot}(X, \mathrm{End}(\V_{\rho}))= H^{\hdot}(E^{\hdot})$  the structure of a Hodge graded Lie algebra described in  example \ref{hdga}.
\end{lem}

\section{Goldman-Millson theory}
In what follows, we fix 
$(\mathbb{V}_{\rho},\mathcal{F}^{\bullet},\overline{ \mathcal {G}}^{\bullet},
  S)$ a $\C$-VHS on $X$ and denote by $\rho:\pi_1(X,x)\to GL(\mathbb{V}_{\rho,x})$ its holonomy.

In \cite{GM} the complete local ring  $\Or$ 
is described rather precisely. We shall review this theory, pointing at some additionnal 
facts easily deduced from the text. We shall also review some additions made by these authors in \cite{GM2}.

\subsection{General representability criteria in Deligne-Goldman-Millson theory}

\subsubsection{The Deligne-Goldman-Millson groupo\"id attached to a dgla}

Given $\mathcal G$ a small groupo\"id, we will denote by $\mathrm{Iso} \ \mathcal{G}$ the set  of its isomorphism classes, i.e.: the quotient  of the set of objects
by the equivalence relation induced by the arrows\footnote{We are working in the category of sets enriched with a functorial version of all small limits and colimits.}. $\mathrm{Iso}$ is a covariant functor from the category of small groupo\"ids
to the category of sets.

Let $L^{\hdot}$  be a dgla, graded by the nonnegative integers and defined over $\C$. We will assume finite dimensionnality of its cohomology objects. Let $(A,\m)$ be a Artin local $\C$-algebra (we assume $A/\m=\C$). Then, one defines the small groupo\"id $DGM (L^{\hdot},A)$ by:

$$
\begin{array}{ccl}
\mathrm{Obj} \  DGM(L^{\hdot},A)& =& \{ \alpha \in L^1\otimes \m | \ d\alpha+ \frac{1}{2}[\alpha,\alpha]=0 \}\\
\mathrm{Hom}_{DGM(L^{\hdot},A)} (\alpha, \beta) & = & \{ \lambda \in L^0\otimes \m | \ \exp(\lambda)\alpha=\beta \}.
\end{array}
$$

This gives a covariant functor $DGM(L^{\hdot},-)$ on the category $\Ar$ of $\C$-Artin local rings with values in the category of small groupo\"ids. Let us describe this functor. In the sequel, the generic object of $\Ar$ will be denoted as $A$, and  $\m$ will stand for the maximal ideal of $A$. The groupo\"id $DGM(L^{\hdot},A)$ is the transformation groupo\"id \footnote{see \cite{GM} p. 52 for this notion.} associated to the set-theoretic action
of the simply connected nilpotent (infinite dimensionnal) Lie group $\exp(L^0\otimes \m)$ on $\mathrm{Obj} \  DGM(L^{\hdot},A)$\footnote{see \cite{GM}, 2.2, p. 53. The action is defined in loc.cit. 1.3 pp. 50-51 using the Baker-Campbell-Hausdorff formula. The differentiable structure on $\exp(L^0\otimes \m)$ will not be used, only its set theoretic group structure.}. 
The transition maps are the obvious ones.

For every action of a group $H$ on a set $S$, we will denote by $[S/H]$ the associated transformation groupo\"id. 
Hence, $$DGM(L^{\hdot},A)= [ \mathrm{Obj} \  DGM(L^{\hdot},A) / \exp(L^0\otimes \m)].$$ 
If $H$ acts on another set $Y$, we define\footnote{see \cite{GM} p. 63.} $[S/H]\bowtie Y := [S\times Y/H]$.

 \begin{rem} In the notation of \cite{Man}, 
$MC_{L^{\hdot}}(A)=\mathrm{Obj} \  DGM(L^{\hdot},A)$ and $Def_{L^{\hdot}}= \mathrm{Iso} \ DGM(L^{\hdot},\_)$.
 \end{rem}

\begin{prop} \label{kur}
 If $H^0(L^{\hdot})=0$ then $\mathrm{Iso} \ DGM(L^{\hdot},\_)$ is prorepresentable. In general, for every splitting $\delta$ of $L^{\hdot}$, a hull in the sense of \cite{Sc} $Kur_{L^{\hdot}}^{\delta}\to Def_{L^{\hdot}}$ can be constructed which is called the formal Kuranishi space.
\end{prop}

\begin{prf}
 \cite{GM2}, see also \cite{Man}. 
\end{prf}

For the reader's convenience, we need to recall the:
\begin{defi}
Let $F$ be a covariant functors of Artin rings. 
Let $T$ be the formal spectrum of a complete local algebra $R$
and $\phi: T\to F$ is a morphism of functors\footnote{We use
Yoneda lemma to justify the abuse of language of using the same notation for $T$ and its functor of points.} given by $\xi\in F(T)$ . $\phi$ is a hull
if and only if it
has the following properties:
\begin{enumerate}
 \item $\phi(\C[\epsilon]/(\epsilon^2))$ is an isomorphism. 
\item For every $B$ an artin ring and $w\in F(B)$ there is 
a morphism $\psi: R\to B$ such that $w= \psi (\xi)$.
\end{enumerate}

\end{defi}

If unicity occurs in the last property then $F$ is prorepresentable by its hull.

\subsubsection{The Deligne-Goldman-Millson groupo\"id attached to an augmented dgla}
Let $\gu$ be a Lie algebra over $\C$ -viewed as a dgla in degree 0- and  $\epsilon:L^{\hdot}\to \gu$ an augmentation of the   dgla $L^{\hdot}$.
Let $(A,\m)$ be a Artin local $\C$-algebra. Then, one defines the small groupo\"id $DGM (L^{\hdot},\epsilon,A)$ by:

$$
\begin{array}{ccl}
\mathrm{Obj} \  DGM(L^{\hdot}, \epsilon, A)& =& \{ (\alpha, e^r) \in L^1\otimes \m \times \mathrm{exp}(\gu \otimes \m) | \ d\alpha+ \frac{1}{2}[\alpha,\alpha]=0 \}\\
\mathrm{Hom} ((\alpha,e^r) , (\beta,e^s)) & = & \{ \lambda \in L^0\otimes \m | \ \exp(\lambda)\alpha=\beta, \ \exp(\epsilon(\lambda)). e^r = e^s\}.
\end{array}
$$
 In terms of the preceding notations, one has:

$$DGM(L^{\hdot}, \epsilon, A)=DGM(L^{\hdot}, A)\bowtie \exp(\gu \otimes \m)_{\epsilon},
$$
the $\epsilon$ subscript meaning that the gauge group acts via $\epsilon$. 

This gives a covariant functor $DGM(L^{\hdot},\epsilon,-)$ on the category $\Ar$ of $\C$-Artin local rings with values in the category of small groupo\"ids.

\begin{prop}\label{rep1}
 If $\epsilon: H^0(L^{\hdot})\to \gu$ is injective then  $DGM(L^{\hdot},\epsilon,\_)$ is a discrete functor in groupo\"ids, $Def_{L^{\hdot},\epsilon}=\mathrm{Iso} \ DGM(L^{\hdot},\epsilon,\_)$ is  prorepresentable. For every splitting $\delta$ of $(L^{\hdot},\epsilon)$ one can construct an explicit  formal scheme and an isomorphism $Kur_{L^{\hdot},\epsilon}^{\delta}\to Def_{L^{\hdot},\epsilon}$. 
\end{prop}
\begin{rem}
$Kur_{L^{\hdot},\epsilon}^{\delta}\to Def_{L^{\hdot},\epsilon}$ is uniquely determinined up to a unique isomorphism.
\end{rem}

\begin{prf} The proof, implicit in \cite{GM}, is an adaptation of \cite{GM2} with an $\epsilon$. It may nevertheless be useful to give an outline.

Suppose $(K^{\cdot}, d , [\, , \, ] )$ is a positively graded dgla, $\hu$ is a Lie algebra, and $\epsilon : K^0\rightarrow \hu$ a Lie algebra map.
Suppose we are given decreasing filtrations $M^{\cdot}K^{\cdot}$ and $M^{\cdot}\gu$ such that $[M^i, M^j]\subset M^{i+j}$ and 
$d,\epsilon :M^i\rightarrow M^i$. Assume the {\em nilpotency condition} $K^{\cdot} = M^1K^{\cdot}$ and $\hu = M^1\hu$, and $M^k=0$ for $k\gg 0$.
In this case the Lie brackets are nilpotent.

A {\em splitting} is a collection of maps denoted $\delta : K^i\rightarrow K^{i-1}$ such that 
$\delta ^2 = 0$, $d=d\delta d$ and $\delta = \delta d \delta$. In this case, we get a decomposition 
$$
K^i = {\rm im}(d)\oplus {\rm im} (\delta ) \oplus (\ker (d)\cap \ker (\delta )).
$$
Indeed an element $u$ can be written as $u= d\delta (u) + \delta d(u) + (1-d\delta -\delta d)(u)$ and any decomposition into 
$u=du_1+\delta u_2 + u_3$ with $du_3=0$, $\delta u_3=0$ has to be of that form. Conversely, suppose we are given a decomposition
$K^i= A^i \oplus B^i \oplus C^i$ such that  $d(K^{i-1})= A^i$ and $\ker (d) = A^i\oplus C^i$.
Then $d: B^{i-1}\stackrel{\cong}{\rightarrow} A^{i}$ and we can define $\delta |_{A^{i}}$ to be the inverse, extended by $0$ on
$B^i$ and $C^i$. The above decomposition associated to $\delta$ is the same as the given $A^i \oplus B^i \oplus C^i$. So,
to give a splitting is the same thing as to give such a decomposition. 

We assume also given a splitting $\delta : \hu \rightarrow H^0(K^{\cdot})$ such that $\epsilon (u) = \epsilon \delta \epsilon (u)$
for any $u\in H^0(K^{\cdot}) = \ker (d:K^0\rightarrow K^1)$. This is equivalent to specifying a subspace $\ker (\delta )\subset \hu$
complementary to $\epsilon (H^0(K^{\cdot}))$. Injectivity of $\epsilon:H^0(K^{\cdot})\rightarrow \hu$ implies then that $\delta (v) = \delta \epsilon \delta (v)$ 
for any $v\in \hu$. 

We assume that our differentials and splittings strictly preserve the filtration $M^{\cdot}$. This is equivalent to requiring the 
direct sum decomposition to be compatible with the filtration in the sense that the filtered vector space $K^i$ is the direct sum of
the subspaces with their induced filtrations (and same for $\hu$). 

With the above notations, recall the {\em gauge fixing procedure} underlying Goldman-Millson's construction \cite{GM2}:

\begin{lem}
\label{gaugefix}
Suppose $\eta \in K^1$ with $d(\eta ) + \frac{1}{2}[\eta , \eta ] = 0$, and $r\in \hu$. Then there is a unique gauge transformation $e^s$ for
$s\in K^0$, taking $(\eta , e^r)$ to a new Maurer Cartan element $(\zeta , e^z)$ such that $\delta (\zeta )=0$ and $\delta (z)=0$. 
\end{lem}
\begin{prf}
The proof is by induction on the filtration $M^{\cdot}$. Suppose we have chosen our gauge transformation to get to 
$(\zeta , e^z)$ with $\delta (\zeta ) \in M^kK^0$ and $\delta (z)\in M^kH^0(K^{\cdot})$. 
Then use the gauge transformation $s= -\delta (\zeta ) - \delta (z)$. Modulo $M^{k+1}$ the effect of this
gauge transformation adds $d(s)=-d\delta (\zeta )$ to $\zeta$, and adds $\epsilon (s)= -\epsilon \delta (z)$ to $z$.
But $\delta (\zeta - d\delta (\zeta ))= 0$ and $\delta (z-\epsilon \delta (z))=0$, so the new element $(\zeta _1,e^{z_1})$
has $\delta (\zeta ) \in M^{k+1}K^0$ and $\delta (z)\in M^{k+1}H^0(K^{\cdot})$. This proves existence of the gauge transformation
in question. Uniqueness is proved similarly: if $(\eta , e^r)$ already satisfies the gauge-fixing property and $e^s$ is a gauge transformation
taking it to another $(\zeta , e^z)$ which satisfies this property, then we can prove by induction that $s\in M^kK^0$ for all $k$, so $s=0$. 
\end{prf}

Let us turn to the construction of the representing formal scheme for the functor $A\mapsto \mathrm{Iso} DGM(L,\epsilon ; A)$. Choose splittings as described previously,
for the dgla $L^{\cdot}$.
In the analytic case, one possible choice of $\delta $ can be obtained by composing $d^{\ast}$ with the Green's operator. 
Note that $\delta$ will not be in any way compatible with the Lie bracket.  However, it induces morphisms also denoted by $\delta$ on 
the complexes $L\otimes \m _A$ for any artinian algebra $A$. We apply the gauge-fixing procedure to the dgla $K^{\cdot}:= 
L^{\cdot}\otimes \m _A$ with filtration $M^kK^{\cdot}:= L^{\cdot}\otimes \m _A^k$. This filtration will be compatible with the splitting
$\delta$ and satisfies the required nilpotency condition. 

Similarly choose a map $\delta : \gu \rightarrow H^0(L)$ which again induces a splitting on $\hu := \gu \otimes \m _A$,
which is compatible with the filtration $M^k\hu := \gu \otimes \m _A^k$. 

By Lemma \ref{gaugefix} applied with $K^{\cdot}=L^{\hdot} \otimes \m$ $\hu=\gu \otimes \m $, 
any Maurer-Cartan element $(\eta , e^r)$ is gauge-equivalent to a unique Maurer Cartan element $(\zeta , e^z)$ such that $\delta (\zeta )=0$ and $\delta (z)=0$
in addition to the Maurer Cartan equation
$d(\zeta ) + \frac{1}{2} [\zeta , \zeta ] = 0$. By hypothesis $(\ast )$ the gauge transformation is unique. In this way we obtain a functor 
$$
A \mapsto DGM^{\delta }(L,\epsilon ; A):= 
$$
$$
\{ (\zeta ,e^r): \;\; r\in \gu \otimes \m _A, \;\; \zeta \in L^1\otimes \m_A, \;\; d(\zeta ) + \frac{1}{2} [\zeta , \zeta ] = 0, \;\; \delta (\zeta )= 0, \delta (r)=0 \} .
$$
This is given by algebraic equations but in a possibly infinite dimensional space $L^1$. We can reduce to a finite-dimensional space in the following way.
Let $H^1:= (\ker (d)\cap \ker (\delta ))\subset L^1$. We can solve the equations
$$
\delta d(\zeta ) + \frac{1}{2} \delta [\zeta , \zeta ] = 0 \;\; \mbox{in}\;\; {\rm im}(\delta : L^2\rightarrow L^1)
$$
because $L^1$ maps surjectively by $\delta d$ to ${\rm im}(\delta )$. The space of solutions maps
isomorphically to $H^1$ by the projection $P^1:L^1\rightarrow H^1$ which vanishes on ${\rm im}(d)+{\rm im}(\delta )$. Hence the functor of formal solutions
in an artinian local algebra $A$ 
is represented by the formal completion of the vector space $H^1$ at the origin. Let $\widehat{H}^1 \times \widehat{\Omega}$ denote the space of solutions
$(\zeta , e^r)$ with also $\delta (r)=0$. Let $P^2: L^2\rightarrow H^2$ be the projection
vanishing on the images of $d$ and $\delta$. Then the map $\zeta \mapsto P^2([\zeta , \zeta ])$ is a formal regular function from $\widehat{H}^1$ to
$H^2$, and $DGM^{\delta}(L,\epsilon )$ is the zero-set of this map. This gives an explicit representation of the functor. 
This finishes the proof of Proposition \ref{rep1}. 
\end{prf}

\subsection{Goldman-Millson isomorphisms}\label{gm}

\subsubsection{Deforming representations of a K\"ahler group}

Let $x\in X$ and  $$\epsilon_x: H^{\hdot}(X, \mathrm{End}(\V_{\rho}))\to \mathrm{End}(\V_{\rho})_x$$ be the augmentation of the Hodge-dgla
$H^{\hdot}=H^{\hdot}(X, \mathrm{End}(\V_{\rho}))$ by evaluation at $x$.
We also define 
for every Artin local algebra $A$, a group
 $G^0_A= \ker(GL(\V_{\rho,x})(A) \to GL(\V_{\rho,x} ))$\footnote{see \cite{GM}, 3.9 p. 66.
$G^0_{-}$ is a prorepresentable group object in the category of covariant set valued functors on $\Ar$.}. 
Then, we have a morphism of groups $\exp(H^0 \otimes \m) \to G^0_A$ induced by the augmentation map. The left action gives a set-theoretic action of the former group on the latter.

\begin{prop} \label{forma1}
The functor on $\Ar$, $h_{GM}$, defined by  $$A\mapsto \mathrm{Iso} \ DGM(H^{\hdot}, A)\bowtie G^0_A$$ is canonically isomorphic 
to the pro-Yoneda functor on $\Ar$ associated to $Spf(\Or)$, the germ at $\rho$ of $\repx$, i.e. to the functor
$h_{\repx_{\rho}}$:

 $$(A\mapsto \repx (A)_{\rho} =\{ \rho_A \in \repx(A)| \ \rho_A \text{mod} \ \mathfrak{m} = \rho \}).$$ 
\end{prop}

\begin{prf}
 This is actually what is proven in \cite{GM}. The combination of theorem 3.5 p. 63, theorem 6.8 p. 82 and the argument made in sections 7-8 imply that the first functor prorepresents $Spf(\Or)$. But the isomorphism of functors constructed in loc. cit. is actually canonically defined. Since we want to be as explicit as possible, we shall give a more detailled sketch of their argument. 

Let $M^{\hdot} \subset E^{\hdot}$ be the subdgla of $D^c$-closed twisted forms. Then the natural projection $M^i\to H^i_{D^c} (X,\mathrm{End}(\V_{\rho}))$ and the Hodge-theoretic canonical isomorphism 
$H^i_{D^c}(X,\mathrm{End}(\V_{\rho}))\to H^i(X,\mathrm{End}(\V_{\rho}))$ gives rise to a dgla morphism from $M^{\hdot}$ to $H^{\hdot}$. Now, the $D'D''$-lemma implies (see \cite{GM}, sect. 7 or \cite{DGMS}) that:
$$ E^{\hdot} \longleftarrow M^{\hdot} \longrightarrow H^{\hdot}
$$
is a dgla quasiisomorphism. 
Using the $DGM(-,-)$ functor and \cite{GM}, cor 2.12 p. 59, we get a  diagram of objectwise 
equivalences of functors from $\Ar$ to $\underline{\mathrm{Gpd}}$:
$$DGM(E^{\hdot}, -) \longleftarrow DGM(M^{\hdot}, -) \longrightarrow DGM(H^{\hdot}, -).
$$
Since the above quasiisomorphism commutes with  the natural augmentation  $\epsilon_x:E^{\hdot}\to  \mathrm{End}(\V_{\rho})_x$, we deduce\footnote{see \cite{GM} sections 3.7-3.9 p.63-64} objectwise 
equivalences:
$$ DGM(E^{\hdot}, -) \bowtie G^0_{-}\longleftarrow DGM(M^{\hdot}, -) \bowtie G^0_{-}\longrightarrow DGM(H^{\hdot}, -)\bowtie G^0_{-}.
$$

Now, an object $\alpha$ of $DGM(E^{\hdot}, A)$ can be interpreted as a flat $A$-linear connections $D_{\alpha}$ on $C^{\infty}(\V_{\rho} \otimes A)$
such that $D_{\alpha}= D_{\V_{\rho}}\otimes \mathrm{id}_A \ mod \ \m$. Using the identification
$\V_{\rho,x}\to \C^N$ that is built in our hypotheses we see that holonomy defines a map from $\mathrm{Obj} \  DGM(E^{\hdot}, A)$ to $\repx (A)_{\rho}$. \cite{GM} sect. 6 contains (implicitely) that this can be enhanced to a natural equivalence of groupo\"ids:
$$ DGM(E^{\hdot}, A) \to [\repx (A)_{\rho}/G^0_A].
$$

Eliminating the superfluous isotropy using identity as a base point for $G^0_A$ gives a natural  equivalence of groupo\"ids:

$$ [\repx (A)_{\rho}/G^0_A] \bowtie G^0_A \longleftarrow  [\repx (A)_{\rho}/\{ \mathrm{Id} \}].
$$

Since we have a natural identification $$\mathrm{Iso} \ [\repx (A)_{\rho}/\{ \mathrm{Id} \}]=\repx (A)_{\rho},$$  passing to $\mathrm{Iso}$ in the above long chain of objectwise natural equivalences gives the required natural isomorphism. 
\end{prf}

\begin{coro}
 Let $\Ar_n$ be the full subcategory of $\Ar$ whose objects $(A, \m)$ satisfy $\m^{n+1}=0$. Then the Yoneda functor
on $\Ar_n$ corepresented by $\Orn$ 
is canonically isomorphic to the restriction to $\Ar_n$ of 
$$A\mapsto \mathrm{Iso} \ DGM(H^{\hdot}, A)\bowtie G^0_A.$$
\end{coro}

\subsubsection{Parameters for Goldman-Millson isomorphisms}

Let $\Omega^0_{H^1_{\rho}}$ be the formal germ at $0$ of 
$\mathrm{End}(\V_{\rho})_x/\epsilon(H^0(X, \mathrm{End}(\V_{\rho})) )$ and $(S_1,\m_1)$ be its complete local algebra. 
Let $Q_{H^1_{\rho}}$ be the formal germ at $0$ of the quadratic cone 
$$\{ \alpha \in H^1 (X, \mathrm{End}(\V_{\rho})) | \ [\alpha,\alpha]=0 \},$$
and $(S_2,\m_2)$ be its complete local algebra. 

Let $h_1$ and $h_2$ be the pro-Yoneda functors of $\Omega^0_{H^1_{\rho}}$
and $Q_{H^1_{\rho}}$.

\begin{defi} A formal subscheme $\mathfrak{t}\subset (\mathrm{End}(\V_{\rho})_x,0)$ transverse   to $\epsilon(H^0)$ will be called a GM-transversal.

\end{defi}

A GM-transversal $\mathfrak{t}$ is smooth or equivalently its  complete local ring is isomorphic to a ring of power series. In fact, we will tacitely assume  that the transversals we consider are attached to a linear subspace, although more general choices could be made if useful. 

\begin{theo} To every  GM-transversal $\mathfrak{t}$, 
we can associate  a well defined isomorphism of covariant functors on $\Ar$:
  $$GM_{\mathfrak t}: h_{\repx_{\rho}}\buildrel{\sim}\over{\to} h_1 \times h_2.$$

With this choice, the inclusion of $h_1$ in $h_{\repx_{\rho}}$ is the transformation of functors induced
by the inclusion of the formal germ at $\rho$ of $\Omega_{\rho}$ into the formal germ $ \repx_{\rho}$. 
\end{theo}

\begin{prf}
In \cite{GM}, lemma 3.10, an isomorphism of functors between $h_{GM}$ and $h_1\times h_2$
is constructed. By definition:

\begin{eqnarray*}
 h_1(A) &=& \mathrm{End}(\V_{\rho})_x\otimes \m /\epsilon(H^0(X, \mathrm{End}(\V_{\rho})) )\otimes \m) \\
h_2(A) &=& \mathrm{Obj} \ DGM(H^{\hdot}, A) \\
h_{GM}(A) &=& h_2(A) \times \exp(\mathrm{End}(\V_{\rho})_x\otimes \m)/\exp(\epsilon(H^0(X, \mathrm{End}(\V_{\rho}))) \otimes \m)
\end{eqnarray*}

Although this point is not made explicitely in \cite{GM}, one way to define precisely the isomorphism  is to prescribe the additional choice of a formal subscheme $\mathfrak{t}\subset (\mathrm{End}(\V_{\rho})_x,0)$ transverse to $\epsilon(H^0)$. 
The exponential gives an isomorphism $h_{\mathfrak{t}}\to h_1$ and a natural isomorphism $i:h_{\mathfrak{t}}\times h_2 \to h_{GM}$ by:
$$ i( t, \eta_2)= (\eta_2, \exp(t)) /\exp(\epsilon(H^0(X, \mathrm{End}(\V_{\rho}))\otimes \m).
$$
 This choice of a Goldman-Millson isomorphism has indeed the required property, by inspection of the construction.
  \end{prf}
  
Since
$\epsilon(H^0)$ is actually a sub-Hodge structure of the polarized Hodge structure $\mathrm{End}(\V_{\rho})_x$ the
germ $\mathfrak{t}_H$ of the orthogonal complement is a canonical choice for $\mathfrak{t}$, we will call this choice the Hodge transversal.
Admittedly, this is a rather artificial way of rigidifying Goldman-Millson's theorem. Hence, it seems abusive to call the resulting isomorphism
canonical and we will call $GM^c_{\mathfrak{t}_H}$  the preferred Goldman-Millson isomorphism.

\begin{defi} \label{can} The isomorphism $GM^c_{\mathfrak{t}_H}$ constructed above $\repx_{\rho} \to \Omega^0_{H^1_{\rho}}\times T$ will be called the preferred
Goldman-Millson isomorphism attached to the $\C$-VHS $\rho$. 
\end{defi}

In a less artificial way,  we have:

\begin{rem}
The functor on $\Ar$ $h_1'$ defined by  $$A\mapsto \exp(\mathrm{End}(\V_{\rho,x})\otimes \m)/ \exp(\epsilon(H^0(X,\mathrm{End}(\V_{\rho}))\otimes \m))$$
is isomorphic to $h_1$ and the projection map  $h_{GM}\to h'_1$ is a retraction of the inclusion $h'_{1}\to h_{GM}$ given 
by $A\mapsto (0,h'_1(A))/ \exp(\epsilon(H^0(X,\mathrm{End}(\V_{\rho}))\otimes \m))$.
\end{rem}

Let us restate this in terms of complete local algebras.

\begin{coro} Choose a GM-transversal $\mathfrak{t}$. 
 There is an isomorphism   $i^c_{\mathfrak t}:\Or \to S_1\hat \otimes S_2$. The ideal $ \jo$ is mapped by this isomorphism to $S_1\hat\otimes \m_2$. Passing to the quotient,  it induces an isomorphism $\Orn \to S_1\hat\otimes S_2 /\m^n$
 which will be also denoted by $i^c_{\mathfrak t}$ .
\end{coro}
\begin{prf}
This follows from the previous theorem,  and the construction of \cite{Sc}.
\end{prf}

We define the ideal $\mathfrak{q}$ as the ideal mapped to $ \m_1\hat\otimes S_2$ by the dual of the preferred
Goldman-Millson isomorphism. $(\jo^n)_{n\in\N}$ induces on
the complete local ring $\Ot=(\Or/\mathfrak{q},\m)$ the filtration $(\m^n)_{n\in\N}$ . The formal subscheme  $T=Spf(\Ot)\subset Spf(\Or)$ defined by $\mathfrak{q}$ will be called the preferred formal slice at $\rho$. $\Ot$ is canonically isomorphic to $S_2$. 

\subsubsection{Interpretation in terms of formal Kuranashi space}

Since the dgla with a trivial differential $(H_{D^c}^{\hdot}(X,\mathrm{End}(\V_{\rho})),0)$ has just one splitting namely $\delta=0$,  $T$ is isomorphic through a uniquely defined tautological isomorphism to the uniquely defined Kuranishi formal space of $H_{D^c}^{\hdot}(X,\mathrm{End}(\V_{\rho}))$. 

Fixing a K\"ahler metric on $X$, Hodge theory enables us to define the orthogonal  splitting of the dgla $E^{\hdot}=E^{\hdot}(X,\mathrm{End}(\V_{\rho}))$. The preceding chain of quasiisomorphisms gives rise to a 
natural map $\kappa^c: T\to Def_{E^{\hdot}}$ which is a hull. 

Holonomy at $x\in X$ gives a canonical isomorphism $h_x: Def_ {E^{\hdot}} \to Def(\rho)$
where $Def(\rho)$ is the Yoneda functor of Artin rings germ at $\rho$ of the groupo\"id  $\repx//GL(\mathbb{V}_{\rho,x})$ attached to the conjugation action. 

Hence $h^c_x: h_x\circ \kappa^c: T \to Def(\rho)$ is a just a hull of the deformation functor $Def(\rho)$. In particular
for every $g\in \mathrm{Aut}(T/Def(\rho))$,  $h_x\circ\kappa \circ g$ is another hull. One easily establishes using \cite{Rim1} (or \cite{Rim2}):

\begin{lem} \label{ambig}
 $\mathrm{Aut}(T/Def(\rho)) \subset \ker (\mathrm{Aut}(\Ot, \m) \to \mathrm{Aut}(\Ot/ \m^2))$ is the subgroup
$\exp(H^0 (X,\mathrm{End}(\V_{\rho}))\otimes \m^2)$. 
\end{lem}

Observe that $H^0 (X,\mathrm{End}(\V_{\rho}))$ is nothing but the Lie algebra of the reductive algebraic group $H\subset GL(\V_{\rho,x})$ which normalizes $\rho$. $H$ acts as automorphisms of the local system $\mathbb{V}_{\rho}$
hence acts on $T$. On the other hand, choosing a vector space supplementary to $\epsilon(H^0(X,\mathrm{End}(\V_{\rho})))$ in 
$\mathfrak{gl}(\V_{\rho,x}))$ (call  $\mathfrak{t}$ the resulting linear transversal) is tantamount to choosing a splitting
of the cohomology augmented Lie algebra. Hence the choice of a linear transversal gives  us a uniquely defined natural isomorphism of punctual formal schemes:
$$  \mathcal{GM}^c: \widehat{GL(\V_{\rho,x})\times_H T}_{\rho}\to\repx_{\rho}=Spf(\Or). 
$$

However, no choice is necessary when  $H^0(X,\mathrm{End}(\V_{\rho}))=0$ which occurs 
iff the normalizer of $\rho $ finite or in the opposite case when  $\rho$ is the trivial representation.

In general, the best basepoint free way to state the main result of \cite{GM} is that there is a canonical equivalence of formal groupo\"ids:

$$ [T/H]\to [\repx/GL(\V_{\rho, x})]_{\rho}
$$

where the second groupo\"id is the germ at $\rho$ of the conjugation groupo\"id. This is actually independant of $x$ (or more precisely is independant on $x$ up to a unique isomorphism). 

\subsection{Split $\C$-MHS on $\Or$ and $\Ot$}\label{mhs}

Since $H^1(X, \mathrm{End}(\V_{\rho}))$ is a weight one $\C$-HS, $H^2(X, \mathrm{End}(\V_{\rho}))$ is a weight one $\C$-HS,  and $$p=[-;-]: S^2H^1(X, \mathrm{End}(\V_{\rho}))\to H^2(X, \mathrm{End}(\V_{\rho}))$$
 respects the Hodge structure, we deduce that $H^1(X, \mathrm{End}(\V_{\rho}))^*$ is a $\C$-HS of weight $-1$ and 
$$I_2= \mathrm{Im}(p^t:H^2(X, \mathrm{End}(\V_{\rho}))^*\to S^2 H^1(X, \mathrm{End}(\V_{\rho}))^*)$$
 is a weight $-2$ $\C$-HS. Hence 
$$I_n= I_2 S^{n-2}H^1(X, \mathrm{End}(\V_{\rho}))^*\subset S^n H^1(X, \mathrm{End}(\V_{\rho}))^*$$ is a weight $-n$ sub $\C$-Hodge structure. 
Hence, $$\Pi_n = S^n H^1(X, \mathrm{End}(\V_{\rho}))^*/I_n \leqno{(\dagger)}$$ is a weight $-n$ $\C$-HS and the algebra $S_2^H:=\sum_{n\in \N} \Pi_n$  viewed as an infinite dimensional split $\C$-MHS with an algebra structure repsecting the MHS. The weight filtration
is exactly $W_{-n} S_2^H= \m_2 ^n$ since $\m_2= \sum_{n>0} \Pi_n$ is a maximal ideal of the complete local algebra $S_2^H$ and $\Pi_1$ generates $\m_2$ . 

The formerly introduced complete local algebra $S_2$ is canonically isomorphic to the algebra obtained from $S_2^H$ by forgetting the Hodge filtrations, hence $\Ot$ carries a canonical split $\C$-MHS whose weight filtration is given by the powers of the maximal ideal. 

On the other hand $S_1$ is the free complete local algebra generated by the weight $0$-HS $(End(\V_{\rho})/ \epsilon(H^0(X, \mathrm{End}(\V_{\rho}))))^*$. It can thus be viewed as a weight zero $\C$-Hodge algebra. 

Hence $(S_1\hat \otimes S_2)$ carries a canonical split $\C$-MHS whose weight filtration is given by the powers of $S_1\hat\otimes m_2$.

\begin{prop} \label{hs}
The filtration $W_{\bullet}$ is the weight filtration of a split $\C$-Mixed Hodge Structures on $\Orn$. Passing to the limit we get a MHS on $\Or$.\end{prop}
 
\begin{prf}
 Choose a $\mathfrak{t}$ a GM-transversal. Using 
the isomorphism $GM_{\mathfrak t}$, we transfer this mixed Hodge structure to $\Or$, the weight filtration 
being defined by the powers of $\jo$. The power $\m^n$ of the maximal ideal is easily seen to correspond to the split sub $\C$-MHS $\sum_{k+l\ge n} S^k H^0(X,\mathrm{End}(\V_{\rho}))\otimes \Pi_l$.  Thus, the powers of $\jo$ are the weight filtration of a $\C$-MHS on $\Orn$.
 
This concludes the proof of Proposition \ref{hs}.
\end{prf}

\section{Mixed Hodge theoretical aspects of Goldman-Millson's theory at a $\C$-VHS}

The following paragraph continues our exegesis of \cite{GM} and aims at making more explicit the use of formality. Despite the fact that the preferred Goldman-Millson isomorphism is rather artificial, $T$ 
does not depend on any choices, so an explicit description should 
be available.

 An explicit construction of three  \lq canonical\rq \  hulls will be described next. They will not   co\"incide in general.  
One of them orresponds to the preferred Goldman-Millson isomorphism constructed above. 

The new two  slices will be used in the last sections  for producing the $\C$-MHS and $\C$-VMHS we are seeking.

\subsection{Preliminary remarks and definitions}

\subsubsection{Universal Maurer-Cartan elements}

Assume $T=Spf(O_T)$ is a formal scheme, $L^{\hdot}$ a dgla and $h: T\to Def_{L^{\hdot}}$ is a hull. In order to define $h$ 
one needs to construct a {\em universal Maurer Cartan element} for this hull, i.e.
an element of $\mathrm{Obj} \ DGM(L^{\hdot}, O_S)$ whose gauge equivalence class gives rise to $h$. 

\begin{defi}
 Assume $T=Spf(O_T)$ is a formal scheme, $L^{\hdot}$ a dgla a  universal Maurer Cartan element for
$L^{\hdot}$ defined on $T$ is an element of $\mathrm{Obj} \ DGM(L^{\hdot}, O_T)$ whose gauge equivalence class gives rise to
a hull. 
\end{defi}

Observe that the hull $Kur^{\delta}_{L^{\hdot}}$ canonically attached to a splitting $\delta$ 
as in Proposition \ref{kur}
carries a tautological universal Maurer Cartan element.
\begin{lem}\label{univmc}
$\eta \in \mathrm{Obj} \  DGM(L^{\hdot},T)$ be universal Maurer Cartan element and 
$\eta'\in \mathrm{Obj} \  DGM(L^{\hdot},S)$ another Maurer Cartan element. Then there exists 
$\phi: T\to S$ a ring morphism and $e^r\in \exp(L^0\otimes \m_S)$ 
a gauge transformation such that $\eta'=e^r. \phi(\eta)$.

Unicity for $(\phi,r)$ holds if $H^0(L^{\hdot})=0$. 
\end{lem}

If $T$ is an object of $\underline{\mathrm{Art}}_n$ (i.e.: satisfies $\m_T^{n+1}=0$) and $(T,\eta)$ satisfies the preceding universal 
property with respect  to $\underline{\mathrm{Art}}_n$ then we say that $(T,\eta)$ is {\em
a  universal Maurer Cartan element of order $n$}.

\subsubsection{Choice of a model}

The construction on $\mathcal{GM}^c$ depends actually on a choice which is the choice of 
$(\ker (D^c), D)$ to set up a quasiisomorphism of $E^{\hdot}(X,\mathrm{End}(\V_{\rho}))$
with its cohomology algebra. This is fine if $\rho$ is actually real for a certain real structure since $\mathcal{GM}^c$ will be defined over the reals. 
The real Zariski closure of the monodromy group being Hodge, it carries a real structure, but it is not clear how the 
real structure extends in non Zariski dense case. If we drop the reality constraint, other very natural
choices are available. 

We may indeed choose $(\ker(D'), D)$ and $(\ker(D''), D)$ as intermediate models resulting
in two alternative formal Kuranishi spaces $h'_x: T\to Def(\rho)$ and $h_x'': T\to Def(\rho)$.

This gives rise to automorphisms of $T$, or,  to be precise, to the following invariants. The {\em fundamental automorphism couple}
of $T$ is the $Aut(h^c_x)$ conjugacy class of $((h^c_x)^{-1}\circ h'_x, (h^c_x)^{-1}\circ h''_x)$.
The {\em fundamental automorphism couple of $\Or$} is $((\mathcal{GM}^c)\circ(\mathcal{GM}')^{-1},\mathcal{GM}^c\circ(\mathcal{GM}'')^{-1})$.

A reality condition will translate in the property that the elements of a fundamental automorphism couple
can be chosen to be conjugate over the reals.

\subsection{Constructing $\C$-VMHS}\label{vmhs}

In this section we first construct a universal Maurer Cartan element for $T$ and $(\ker(D'),D)$
then for $(\ker(D''),D)$ and glue them together. 

\subsubsection{The $D'D''$-lemma}
\begin{lem} \label{a1} Let $\eta_1, \ldots , \eta_b \in E^{\hdot}$ form a basis of the subspace $\mathcal{H}^1(X,\mathrm{End}(\V_{\rho}))$ of harmonic 
twisted one forms, each $\eta_i$ being of pure Hodge type $(P_i,Q_i)$ for the Deligne-Zucker
$\C$-CHM structure on $E^{\hdot}$. Then  $\{ \eta_i \}$ is a basis of
$H^{1}(X,End(\V_{\rho}))$ whose dual basis we denote by $(\{\eta_1\}^*, \ldots, \{\eta_b^*\})$. 

The $End(\V_{\rho})\otimes\Pi_1 $-valued one-form   $\alpha_1$ defined by \footnote{Recall that $\Pi_1=H^{1}(X,End(\V_{\rho}))^*$ by $(\dagger)$.}:

$$\alpha_1=\sum_{i=1}^b \eta_i \otimes \{\eta_i\}^*
$$

is $D$ and $D^c$-closed. 
 \end{lem}

\begin{prf} This is a consequence of the already mentioned Deligne-Zucker construction \cite{Zuc}.
\end{prf}

Let $\pi_2: \Pi_1 \otimes \Pi_1 \to \Pi_2$ the product mapping constructed above (see $[\dagger]$)

\begin{lem}
 $ \beta_2=\pi_2 ([\alpha_1,\alpha_1]) $ is $D$-closed, $D^c$ closed and cohomologous to $0$.
\end{lem}
\begin{prf}
$[\alpha_1,\alpha_1]\in E^2 \otimes \Pi_1 \otimes \Pi_1$ is $D'$ and $D'$ closed since $E^{\hdot}$ is a dgla and $D'$ also preserves the Lie bracket. 
$\pi_2$ being a linear map, $\beta_2$ is $D'$ and $D''$-closed. 

Let $\phi:\Pi_2 \to \C$ be a linear form. Then $\phi$ can be interpreted as an element $\sum_{1\le i,j\le b} c_{i,j} \{\eta_i\}\otimes \{\eta_j\}$ in the kernel of $$[-,-]: \otimes^2H^1(X,End(\V_{\rho}))\to H^2(X,End(\V_{\rho})),$$ ie. such that:

$$ \sum_{1\le i,j\le b} c_{i,j} [\{\eta_i\}, \{\eta_j\}]=0. \leqno{(R)},
$$

Hence $\phi(\beta_2)= \sum_{1\le i,j\le b} c_{i,j} \eta_i\wedge\eta_j$.
If $(c_{i,j})$ is any antisymmetric matrix, then this form vanishes and $(R)$ holds. If $(c_{i,j})$ is symmetric and satisfies $(R)$ then:
$$\phi(\beta_2)=\frac{1}{2}\sum_{i,j} c_{i,j}[\eta_i,\eta_j].$$

Since $[\eta_i,\eta_j]$ is a De Rham representative of $\{\eta_i\},\{\eta_j\}$
$\phi(\beta_2)$ is cohomologous to zero thanks to $(R)$. Hence $\beta_2$ is cohomologous to zero.

\end{prf}

\begin{lem} There is a form $\gamma_2 \in E^{0}\otimes \Pi_2$ such that 
$D'D'' \gamma_2= \beta_2$. We define $\alpha_2 \in E^1 \otimes \Pi_2$
to be $\frac{1}{2}D'\gamma_2$. $\alpha_2$ is $D'$-exact and satisfies $D\alpha_2+\frac{1}{2}\beta_2=0$.
 \end{lem}
\begin{prf}
This is a consequence of the $D'D''$-lemma, see for instance \cite{DGMS}.
\end{prf}

If $\alpha\in E^1\otimes \Pi_b$, $\beta\in E^1\otimes \Pi_b$ we define $\alpha\wedge \beta \in E^2 \otimes \Pi_{a+b}$ using composition of matrices and 
product in $\Pi$. With this notation $\alpha_1\wedge \alpha_1 = \frac{1}{2} \beta_2$. 

\begin{lem} $\beta_3=\alpha_1 \wedge \alpha_2+\alpha_2\alpha_1$ is $D'$-exact and $D''$-closed. By the $D'D''$ lemma we conclude $\beta_3= D'D''\gamma_3$ and we define 
$\alpha_3=D'\gamma_3$. 
 \end{lem}
\begin{prf}
$D\alpha_2+\alpha_1\alpha_1=0$. Hence
 $D\beta_3=\alpha_1\wedge \alpha_1 \wedge \alpha_1-\alpha_1\wedge \alpha_1 \wedge \alpha_1$. 
\end{prf}

\begin{lem} 
 For $k\ge 3$ we can construct a $D'$ exact form  $\alpha_k \in E^1\otimes \Pi_k$  such that the following relation holds:
$$
D\alpha_k +\alpha_{k-1}\alpha_1 + \alpha_{k-2}\alpha_2 +\ldots+ \alpha_1 \alpha_{k-1}=0
$$
\end{lem}

\begin{prf}
 We have done the $k=3$ case explicitely. The first step for $k\to k+1$ is to 
observe that $\alpha_k\alpha_1 + \ldots +\alpha_1 \alpha_k$ is $D''$ closed and $D'$-exact and then apply the $D'D''$ lemma. 
\end{prf}

\begin{lem}
 $\sum_{k=1}^{+\infty} \alpha_k$ is a universal Maurer Cartan element for $T$ and $(\ker(D'),D)$.
\end{lem}
\begin{prf}
 This follows from the form of $\alpha_1$. 
\end{prf}

\begin{rem}
 Since $\mathcal{H}^1(X,\mathrm{End}(\V_{\rho}))=\ker(D':E^1\to E^2)\cap \ker (D'':E^1\to E^2)$ and $\ker(D'D'': E^0 \to E^2)=\ker (D'')$  the construction of $(\alpha_k)_{k\ge 1}$
is  canonical, purely complex analytic (does not depend on the K\"ahler structure)  functorial and the $\alpha_k$ are uniquely determined. 
\end{rem}

\subsubsection{A canonical connection}

Let $(V_n, D)$ the smooth vector bundle with a flat connection 
underlying the local system :
$$\V_{\rho}\otimes_{\C} \Otn=\oplus_{k=-n}^0 \V_{\rho} \otimes \Pi_{-k}.$$ 
The weight filtration on $\Otn$ gives rise to a filtration  $\{ (W_k V_n,D)\}_{-n \le k\le 0}$. Actually $\V_{\rho}\otimes_{\C} \Otn$ is naturally a split 
VMHS.

Let $A'_k: V_n\to V_n\otimes \mathcal{E}^1_X$ defined as $\sum_l A_k^l$ where  
$A_k^l: \V_{\rho} \otimes \Pi_l \to \V_{\rho} \otimes \Pi_{l+k} \otimes \mathcal{E}^1_X$ is the natural multiplication by $\alpha_k$. 

\begin{prop} \label{expl}
 $D+\sum_{k=1}^n A'_k$ is a flat $\Otn$-linear connection on the filtered smooth vector bundle
 $\{ (W_k V_n)\}_{-n \le k\le 0}$. Its holonomy representation 
$\rho_{T,n}:\pi_1(X,x) \to GL (\V_{\rho,x}\otimes\Otn)$ passes to the limits and gives rise
to $\rho_T:\pi_1(X,x)\to GL(\V_{\rho,x}\otimes\Otn)$ which induces the Goldmann-Millson slice
$T\to \repx$.
\end{prop}

\begin{prf} The fact that the connection is flat is a restatement of the construction in the previous subsection. The link with 
the construction in \cite{GM} can be easily made by inspecting this reference.
\end{prf}

\subsubsection{Griffiths transversality}

Since $ \V_{\rho} \otimes \Pi_{-k}$ carries a weight $-k$ $\C$-HS, $(V_n,D)$
underlies a split VMHS whose Hodge and anti-Hodge bundles are:

$$\mathcal{F}^p(V_n)=\oplus_{k=-n}^0 \mathcal{F}^p(\V_{\rho} \otimes \Pi_{-k})$$

$$\overline{\mathcal{G}}^q (V_n)=\oplus_{k=-n}^0 \overline{\mathcal{G}}^q( \V_{\rho} \otimes \Pi_{-k}).$$

\begin{lem}\label{antitrans}
The connection $D+A'=D+\sum_{k=1}^n A'_k$ is Griffiths  transversal in the sense that
$$ (D+\sum_{k=1}^n A'_k)^{1,0} C^{\infty}( {\mathcal{F}}^q )\subset  C^{\infty}( {\mathcal{F}}^{q-1} )\otimes \Omega^{1,0} $$
$$ 
 (D+\sum_{k=1}^n A'_k)^{0,1} C^{\infty}( {\mathcal{F}}^q )\subset  C^{\infty}( {\mathcal{F}}^{q} )\otimes \Omega^{0,1}.
$$

\end{lem}
\begin{prf}
 Observe that the twisted one-form $\alpha_1$ is a Hodge type $(0,0)$ vector
in a weight $0$ HS.
Since it is in $\mathcal{G}^{0} $ it follows that $A'_1 .{\mathcal{F}}^q (V_n)\subset {\mathcal{F}}^q (V_n\otimes E^1)$, ie. : 
$$A'_1 .{\mathcal{F}}^q (V_n)\subset {\mathcal{F}}^q \otimes\Omega^{0,1} +  {\mathcal{F}}^{q-1} \otimes\Omega^{1,0}.$$

Now, $\beta_2$ is also a Hodge type $(0,0)$ vector
in a weight $0$ HS. Hence, $\gamma_2$ is also a Hodge type $(-1,-1)$ vector
in a weight $-2$ HS. Hence $\alpha_2$ is of Hodge type $(0,-1)$ 
in the weight $-1$ HS $E^1\otimes \Pi_2$. Hence $\alpha_2\in \mathcal{F}^{0}\cap\overline{\mathcal{G}}^{-1} (E^1 \otimes \Pi_2)$. Since it is in $\mathcal{F}^{0} $ it follows that $A'_2 . {\mathcal{F}}^q (V_n)\subset {\mathcal{F}}^q (V_n\otimes E^1)$. 

Continuing this way, we see that $\alpha_k$ is of Hodge type $(0,1-k)$ 
and $$(\sum A'_k).{\mathcal{F}}^q (V_n)\subset {\mathcal{F}}^q (V_n\otimes E^1).$$ 
This is Griffiths transversality. 
\end{prf}

\subsubsection{The second filtration}

\begin{rem}
 $D+(\sum_{k=1}^n A'_k)$ is not Griffiths anti-transversal for the above definition of $\overline{\mathcal{G}}^{\hdot}$.
\end{rem}

So, we need to introduce the following variant of our basic construction:

\begin{lem}\label{trans} Let $\alpha_1^v$ denote the twisted 1-form $\alpha_1$ called by another name. 
 For $k\ge 2$, we can construct a $D''$-exact form  $\alpha^v_k \in E^1\otimes \Pi_k$  such that the following relation holds:
$$
D'\alpha^v_k +\alpha^v_{k-1}\alpha^v_1 + \alpha^v_{k-2}\alpha^v_2 +\ldots+ \alpha^v_1 \alpha^v_{k-1}=0.
$$
$\alpha^v_k$ is of Hodge type $(1-k,0)$ and $D+A''=D+\sum \alpha^v_k $ satisfies Griffiths antitransversality for the above $\overline{\mathcal{G}}^{\hdot}$.
\end{lem}

In order to get a VMHS structure on $\rho_{T,n}$ one needs to prove the following proposition where $\m$ stands for the maximal ideal in $\Otn$:
\begin{prop} \label{gauge}
There is $W_{\hdot}$-preserving gauge transformation $g$  and $\phi^* \in \mathrm{Aut}(\Otn)$ 
such that: 
\begin{enumerate}
 \item $g\in\exp (E^0 \otimes \m^2)$ hence induces $\mathrm{Id}$ on $Gr_W^{\hdot} V_n$,
\item $(\phi^* :\m/\m^2 \to \m/ \m^2)=\mathrm{id}_{\m/\m^2}$, hence $\mathrm{id}_{\V_{\rho}}\otimes \phi^*$ preserves 
$W_{\hdot}$ and induces $\mathrm{Id}$ on $Gr_W^{\hdot} V_n$,
\item $$D+A''= g(\mathrm{id}_{\V_{\rho}}\otimes \phi^*) (D+A') (\mathrm{id}_{\V_{\rho}}\otimes \phi^*)^{-1}g^{-1}.$$
\end{enumerate}
\end{prop}
\begin{prf}

Let us remark first that the automorphisms of $\Otn$ act on the group of gauge transformations $\exp(E^0\otimes \m)$
so that the group of filtered bundle automorphisms they generate is isomorphic to the semi direct product deduced from this action. 

Since $A'_1=A_1$,  the statement is obviously true for $n=1$. 
Since $\alpha_2^v= + \frac{1}{2} D'\gamma_2$, we have  $\alpha_2^v -\alpha_2 = D\gamma_2$. $\exp(\pm\gamma_2)$ is the required gauge transformation for $n=2$, in which case we can still set $\phi^*=\mathrm{id}$.
For larger $n$, we have not been able to get such an explicit construction of the pair $(g,\phi^*)$.

In general, since $A'$ and $A''$ come from a universal Maurer Cartan element, this is a consequence of the universality property given in lemma \ref{univmc}.
\end{prf}

\begin{theo}

 The filtered vector bundle $(V_n, W_{\hdot})$ with connection $D+A'$ constructed in prop \ref{expl}, Hodge filtration constructed in lemma \ref{antitrans}, and anti-Hodge filtration defined by transporting the anti-Hodge filtration
constructed in lemma \ref{antitrans} using the bundle automorphism
$g(\mathrm{id}_{\V_{\rho}}\otimes \phi^*)$ constructed in proposition \ref{gauge} gives rise 
to a $\C$-VMHS whose holonomy is the filtered representation $\rho_{T,n}$.
\end{theo}
\begin{prf}
Since the bundle automorphism contsructed in proposition \ref{gauge} induces identity on $Gr_W^{\hdot}(V_n)$
the new Hodge filtration still defines on each stalk a $\C$-VMHS thanks to lemma \ref{twist}. Griffiths antitransversality  is obtained by transport of structure from lemma \ref{trans}.
\end{prf}

This concludes the proof of Theorem \ref{gln}.

\begin{rem}
 In the real case,  if one insists in using the definition that $\overline{\mathcal{F}^{\hdot}}$ is the complex conjugate
of ${\mathcal{F}^{\hdot}}$,
one has to transport ${\mathcal{F}^{\hdot}}$ to a real model, for instance the model given by $(\ker D^c, D)$
using a real fundamental automorphism couple and a pair of conjugate gauge transformations. 
\end{rem}

\begin{rem}
 In case $H^0(X,End(\V_{\rho}))\not= \C \mathrm{Id}$, this VMHS is NOT uniquely defined, since we may twist it by the action of 
the ambiguity group. 
\end{rem}

\subsection{$\C$-MHS on $\Or$}

We may use the preceding construction and lemma \ref{twist} to put a $\C$-MHS on $\Ot$ which will be only defined up to the action of the ambiguity group $\mathrm{Aut}(T/Def(\rho))$ constructed in Lemma \ref{ambig}. 

Since the reductive group $H$ acts as a group of automorphisms of $\V_{\rho}$ the whole construction is $H$-invariant \cite{Rim2} which yields a uniquely defined
$C$-MHS on $O_T^H$ which is the ring of formal series at $[\rho]$ in $M_B(X,GL(\V_{\rho,x})$. 

To treat the case of $\Orn$, we need to introduce an augmented
version of the preceding construction.

\subsubsection{Filtered Goldman-Millson theory}

Consider the following situation: we are given a nonnegatively graded dgla $L$ with a decreasing filtration $G^{\cdot}$ such that
$d:G^pL\rightarrow G^pL$ and $[\, , \, ] : G^pL\times G^qL \rightarrow G^{p+q}L$. Let $Gr_G(L)= \bigoplus G^pL/G^{p+1}L$ 
be the associated-graded dgla. Suppose furthermore that we are given a finite dimensional Lie algebra
$\gu$ and an augmentation $\epsilon : L^0\rightarrow \gu$ compatible with the Lie bracket. Suppose $\gu$ is also given a filtration
denoted $G^{\cdot}{\gu}$, compatible with Lie bracket, and $\epsilon$ is compatible with the filtrations.

Suppose $B$ is an artinian local algebra also provided with a decreasing filtration denoted $G^{\cdot}B$
compatible with the algebra structure $G^pB\times G^qA\rightarrow G^{p+q}A$, with $1_B\in G^0B$. 
We assume that the filtration is {\em exhaustive}, that is $G^pB=B$ for $p\ll 0$ and $G^pB = 0$ for $p\gg 0$. Let 
$\m _B$ denote the maximal ideal of $B$, which has its induced filtration $G^{p}\m _B := G^pB \cap \m _B$. 

Let $\exp (\gu \otimes \m_B)$ denote the nilpotent Lie group associated to the nilpotent Lie algebra 
$\gu \otimes \m_B$. 

An {\em augmented Maurer-Cartan element} is a pair $(\eta , e^r)$ with 
$$
\eta \in L^1\otimes _{\C}\m _B  , \;\;\; d(\eta ) + \frac{1}{2}[\eta , \eta ]=0,
$$
and $e^r \in \exp (\gu \otimes \m_B)$. Let $\mathrm{Obj} \  DGM(L,\epsilon ; B)$ denote the space of agumented Maurer-Cartan elements. 
An element $(\eta , e^r )$ is {\em compatible with the filtrations} if $\eta \in G^0(L^1\otimes \m_B)$ and
if $\beta \in \exp (G^0(\gu \otimes \m_B))$. Let $\mathrm{Obj} \  DGM(L,\epsilon ; B)^G$ denote the subspace of elements compatible with the
filtration. 

The group $\exp (L^0\otimes \m_B)$, with elements denoted $e^s$, acts as a group of gauge transformations which acts on both components. 

The group $\exp (G^0(L^0\otimes \m_B))$ acts on $\mathrm{Obj} \  DGM(L,\epsilon ; B)^G$. Let 
$$
DGM(L,\epsilon ; B)^G:= [\mathrm{Obj} \  DGM(L,\epsilon ; B)^G/ \exp (G^0(L^0\otimes \m_B))]
$$
denote the {\em filtered Deligne-Goldman-Millson groupoid},  quotient groupoid of the filtered-compatible space
by the filtered gauge group.

Suppose $K$ is another dgla, with filtration $G^{\cdot}K$ and an augmentation $\varepsilon$ towards the same $\gu$.
A {\em filtered augmented quasiisomorphism} from $K$ to $L$ is a morphism of dgla's $\psi : K\rightarrow L$, 
compatible with the filtrations, making a commutative square with the augmentations, and such that 
$Gr_G(\psi ) : Gr_G(K)\rightarrow Gr_G(L)$ is a quasiisomorphism. 

The quasiisomorphism invariance of the Deligne-Goldman-Millson groupo\"id
 generalizes here:

\begin{prop}
\label{invariance}
Suppose $\psi$ is a filtered augmented quasiisomorphism, then the induced map
$$
DGM(K,\epsilon ; B)^G\stackrel{DGM(\psi ; 1_B)}{\rightarrow} 
DGM(L,\epsilon ; B)^G
$$
is an equivalence of groupoids.
\end{prop}
\begin{prf}
Let $M^{\cdot}$ be the filtration of $\m _B$ by powers of $\m _B$. We can 
choose a common splitting 
$\m _B = \bigoplus V^{p,q}$ for the two filtrations $M^{\cdot}$ and $G^{\cdot}$ 
on the vector space $\m _B$,
not necessarily compatible with the algebra structure. Thus $M^k = \bigoplus _
{p\geq k}V^{p,q}$ and 
$G^r= \bigoplus _ {q\geq r}V^{p,q}$. The filtration $M^{\cdot}$ induces 
filtrations going by the same name 
on $K^{\cdot}\otimes \m _B$ and $L^{\cdot}\otimes \m _B$. Furthermore, we can 
express the filtrations $G^{\cdot}$
on these complexes as
$$
G^r(L^{\cdot}\otimes \m _B) = \bigoplus _{p, j} G^j(L^{\cdot})\otimes V^{p,r-j},
$$
$$
G^r(K^{\cdot}\otimes \m _B) = \bigoplus _{p, j} G^j(K^{\cdot})\otimes V^{p,r-j}.
$$
These expressions also give splittings for the filtrations $M^{\cdot}$. 
These expressions are compatible with the differential, but not with the 
bracket. 
However, the condition for a map to be a filtered quasiisomorphism depends only 
on the differential.
Thus, the morphism $K^{\cdot}\otimes \m _B\rightarrow L^{\cdot}\otimes \m _B$ is
a bifiltered quasiisomorphism with respect to the pair of filtrations $G^
{\cdot}, M^{\cdot}$.
In particular, the morphism 
$$
G^0(K^{\cdot}\otimes \m _B) \rightarrow G^0(L^{\cdot}\otimes \m _B)
$$
is a filtered quasiisomorphism for the filtration $M^{\cdot}$, and this 
filtration makes the bracket nilpotent
(that is, the bracket is trivial on the associated graded of $M^{\cdot}$).
A similar discussion holds for the filtrations on $\gu \otimes \m_B$. 

The groupo\"id in the filtered case $DGM(K,\epsilon ; B)^G$ is just the groupoid 
of Maurer-Cartan elements
in $G^0(K^{\cdot}\otimes \m _B)$ together with a framing in $G^0(\gu \otimes \m 
_B)$, and the same for $L^{\cdot}$.
Using \cite{GM}, Cor 2.12, p. 59, we conclude that
$$DGM(K,\epsilon ; B)^G\stackrel{DGM(\psi ; 1_B)}{\rightarrow} 
DGM(L,\epsilon ; B)^G$$ is an equivalence of groupoids. 
\end{prf}

On the other hand, if $f:B\rightarrow B'$ is a morphism of filtered artinian local algebras, we obtain a morphism
$$
DGM(L,\epsilon ; B)^G \stackrel{DGM(1_L; f)}{\rightarrow}
DGM(L,\epsilon ; B')^G. 
$$
Say that an element $(\eta , e^r)\in \mathrm{Obj} \  DGM(L,\epsilon ; B)^G$ is an {\em order $k$ universal augmented filtered Maurer-Cartan element}
if $\m _B^k=0$ and, for any artinian local algebra with exhaustive filtration $(R,G^{\cdot})$ such that $\m _{R}^k=0$ 
the map $f\mapsto DGM(1_L; f)(\eta , e^r)$ induces an equivalence from the discrete groupoid of filtered algebra morphisms
$B\rightarrow R$, to the filtered DGM groupoid 
$$
Hom _ {\rm filt.alg}((B,G^{\cdot}),(R,G^{\cdot}))\stackrel{\sim}{\rightarrow} DGM(L,\epsilon ; R)^G.
$$
Note in particular that this condition means that for any $R$ the objects in the groupoid $DGM(L,\epsilon ; R)$
don't have nontrivial automorphisms.

Let us now develop an filtered analog of Proposition \ref{rep1}.
We don't consider the general question of representability of the filtered DGM groupoid functor. 
Probably, in the very general case representability will not hold, so some conditions on the filtrations would be necessary. 
In the case of interest to us,
using filtered augmented quasiisomorphisms we can reduce to a case where the filtration is decomposed, in which case
it is easier to show representability. 

We say that a filtered augmented dgla $(L, \epsilon , G^{\cdot})$ is {\em decomposed} if there exists an action of $\C^{\ast}$ on
the dgla $L$, and an action on $\gu$ such that the augmentation is compatible with the action, and such that
$L$ and $\gu$ decompose into eigenspaces for the action which split the filtrations $G^{\cdot}$. 
Concretely this means that we are given isomorphisms $L\cong Gr_ G(L)$ and $\gu \cong Gr_G(\gu )$ which are
compatible with the differential, the bracket, and the augmentation. Denote by $L(k)$ the subspace corresponding to $Gr^k_G(L)$
and similarly $\gu (k)$. 
The decomposition corresponds to the action of  $\C ^{\ast}$ on $(L,\gu ,\epsilon )$. 

We can use the decomposition condition to show representability
of the filtered DGM groupoid functor. 

Assume $\epsilon: H^0(L^{\hdot})\to \gu$ is injective. 
Choose a splitting $\delta$ compatible with the decompositions of $L^{\cdot}$ and $\gu ^{\cdot}$. 
We get an action of $\C^{\ast}$ on the representing formal 
scheme $DGM^{\delta}(L,\epsilon )$ constructed explicitly above. Let $R$ be the complete local coordinate ring of $Kur_{L^{\hdot},\epsilon}^{\delta}$,
with its universal Maurer Cartan element $\eta$.  From the construction, we see that the action on $R/\m _{R}^k$, is an algebraic action. Furthermore, the universal Maurer Cartan element
$\eta$ is compatible with these actions of $\C^{\ast}$. 

Let $G^{\cdot}$ denote the filtration of $L^{\cdot}$ and $\gu$ corresponding to the decomposition.
Note that $d$ and $\epsilon$ are strictly compatible with $G^{\cdot}$. The decomposition of $R$ induces a filtration which we also denote by $G^{\cdot}$,
and the universal Maurer Cartan element lies in $DGM(L,\epsilon ; R)^G$. 

\begin{theo}
\label{filtrep}
Suppose $B$ is an artinian local algebra with action of $\C^{\ast}$, and suppose $\zeta \in DGM(L,\epsilon ; B)$ is an Maurer Cartan element fixed by $\C^{\ast}$.
Then the corresponding map $\nu : R\rightarrow B$ and gauge transformation $w$ between $\nu (\eta )$ and $\zeta$, are fixed by $\C^{\ast}$.

Suppose $B$ is provided with a decreasing filtration $G^{\cdot}$ and $(\eta , e^r)\in DGM(L,\epsilon ; B)^G$ is an Maurer Cartan element compatible with the
filtration. Then the map $\nu : R\rightarrow B$ and gauge transformation $w$ are compatible with the filtrations. 
\end{theo}
\begin{prf}
If $B$ has a decomposition and $(\eta , e^r)$ is an Maurer Cartan element preserved by $\C^{\ast}$, then let $e^s$ be the gauge transformation going from $(\eta , e^r)$ to
$(\zeta , e^z)$ with $\delta (\zeta )= 0$ and $\delta (z)=0$. By unicity of $s$, we have that $s$, $\zeta$ and $z$ are fixed by $\C^{\ast}$. The $(\eta ,e^z)$ give
the coordinates for the map $\nu : R\rightarrow B$ and $s$ gives $w$, so $\nu$ and $w$ are fixed by $\C ^{\ast}$.  

Consider now the filtered case. The dgla at index $0$ in the filtration can be expressed as
$$
G^0(L^{\cdot}\otimes \m  _B) = \bigoplus _j L^{\cdot}(-j)\otimes G^j\m _B, \;\;\; G^0(\gu \otimes \m _B) = \bigoplus _j \gu (-j)\otimes G^j\m _B.
$$

The splitting $\delta$ is defined separately on each complex $L^{\cdot}(-j)$ and $\gu (-j)$, so it induces
a splitting of the dgla $K^{\cdot} := G^0(L^{\cdot}\otimes \m _B)$ with Lie algebra
$\hu := G^0 (\gu \otimes \m _B)$. Hence the gauge fixing lemma 2.4 can be applied to $(K^{\cdot}, \hu )$. 

Given an Maurer Cartan element $(\eta ,e^r)\in DGM(L,\epsilon ; B)^G$, this means exactly that we have an Maurer Cartan element for $(K^{\cdot}, \hu )$
so by Lemma \ref{gaugefix} there is a unique gauge transformation $s\in K^0=G^0(L^{0}\otimes \m  _B)$
transforming $(\eta ,e^r)$ to an Maurer Cartan element $(\zeta ,e^z)$ for $(K^{\cdot}, \hu )$ with $\delta (\zeta )=0$ and $\delta (z)=0$.
This new element is again in $DGM(L,\epsilon ; B)^G$, and it corresponds to a morphism $\nu : R\rightarrow B$ which sends 
$R(j)$ to $G^j(B)$.
\end{prf}

\begin{rem}
 When finishing this paper, we realized that the $\C^{\ast}$-invariant part of this construction was actually done in \cite{GK}.
\end{rem}

\begin{coro} \label{filtrep2}
Suppose $(L,\epsilon , G^{\cdot})$ is a filtered augmented dgla, with finite dimensional cohomology groups $H^i$ for 
$i=0,1,2$. Suppose that $\epsilon : H^0L\rightarrow \gu$ is 
injective. 
Suppose that $(L,\epsilon , G^{\cdot})$ is filtered quasiisomorphic to a filtered augmented dgla which has a splitting
of the filtration compatible with differential, bracket and augmentation. 

Then for any $k$ there exists a filtered artinian algebra $R$ with $\m _R^k=0$
and a $k$-th order universal augmented Maurer-Cartan element compatible with filtration $(\eta , e^r)\in \mathrm{Obj} \  DGM(L,\epsilon ; R)^G$. In fact, $R$ may also be split in
the sense that there is an action of $\C^{\ast}$ (or equivalently $R\cong Gr(B)$) and $(\eta , e^r)$ is compatible
with the splitting; this splitting depends on a choice of filtered quasiisomorphism with a split dgla. 

The universal object $(R,G^{\cdot}; (\eta , e^r))$ is unique up to unique isomorphism and gauge transformation. 
That is to say that if $(R',G^{\cdot};(\eta ',e^{r'}))$ is another universal element then there is an isomorphism 
of filtered algebras $\nu : (R, G^{\cdot}) \cong (R', G^{\cdot}) $
and a gauge transformation $s\in \exp (G^0(L^0\otimes \m _{B'}))$ such that $e^s(\nu (\eta )) = \eta '$.  
The pair $(\nu , s)$ such that $e^s(\nu (\eta )) = \eta '$ is
unique.
\end{coro}
\begin{prf}
For the split dgla, Theorem \ref{filtrep} provides the representability. By the 
invariance statement of Proposition \ref{invariance}, this representing object works for the original $L^{\cdot}$,
and unicity comes from the universal property. 
\end{prf}

\begin{coro}\label{univ}
Let $R'$ be an artinian algebra and a $k$-th order universal
Maurer Cartan element $(\eta , e^r)\in DGM(L^{\cdot}, \epsilon ; R)$. Let  $(R, G^{\cdot})$ be as in corollary \ref{filtrep2} and $(\eta ', e^{r'})$ a filtered $k$-th order universal
Maurer Cartan element. 

Then there is a unique isomorphism $\nu : R'\cong R$ 
and a unique gauge transformation $e^w$ going from $\nu (\eta ',e^{r'})$ to $(\eta , e^{r})$. The isomorphism $\nu$ induces a filtration
$G^{\cdot}R' := \nu (G^{\cdot}R)$ on $R'$. 
\end{coro}
\begin{prf}
 Indeed the splitting used to construct the universal filtered Maurer Cartan-element is a splitting that can be used to construct a universal Maurer Cartan-element - the construction being parallel. Hence the filtered Maurer Cartan-element is good enough to serve as an ordinary Maurer Cartan-element.  
\end{prf}

However, the original universal element $(\eta , e^r)$ is not necessarily compatible with the filtration;
we have existence of a gauge transformation $e^w$ such that $e^w\cdot (\eta , e^r)$ is compatible with the filtration.

\subsubsection{The mixed Hodge structure on the formal completion of the representation space}

We now wish apply Proposition \ref{mhalg} to the case $A= \Orn$. Let us first define 3 filtrations on this ring:

\begin{defi} 
 The weight filtration on $\Orn$ is by powers of the ideal $W_{-k}= \jo ^k$.
The Hodge filtration $F$ comes from corollary \ref{univ} applied to the Hodge filtration
 of the augmented  Goldman-Millson DGLA $(E^{\hdot}(X,End(\V_{\rho})), \epsilon_x)$ where the augmentation is evaluation at $x$. The anti Hodge filtration comes from the same construction.
\end{defi}

\begin{theo}
These provide $A=\Orn$ with a CMHS which is unique up to a unique isomorphism.
\end{theo}
\begin{prf}
Let us check the conditions in Proposition \ref{mhalg}. 

Condition (1) comes from Goldman-Millson's theorem (indeed the product of a quadratic
cone with a vector space is again a quadratic cone). 

For condition (2), note that $V$ is the dual of the space of deformations 
of $\rho$ in $R(X,x,G)$. This space of deformations is a relative cohomology group:
$V^{\ast}= H^1((X,x),Ad(\rho ))$. This has a mixed Hodge structure, which is exactly the one given by the restrictions
of the filtrations $W,F,G$ above. There are only two weight quotients, in degrees $0$ and $-1$ for $V$ or
degrees $0$ and $1$ for $V^{\ast}$. 

For condition (3),  note that the kernel is given by the obstruction map. Since $Ad(\rho )$ has a Lie algebra structure
which is antisymmetric, the multiplication $H^1\times H^1\rightarrow H^2$ is symmetric, that is
$$
Sym^2(V^{\ast})=Sym^2(H^1(X,x,Ad(\rho ))) \rightarrow H^2(X,Ad(\rho )).
$$
We didn't include the basepoint $x$ on the right because it doesn't affect $H^2$.
The transpose or dual of this map is 
$$
H^2(X,Ad(\rho ))^{\ast} \rightarrow Sym^2(V)
$$
and by Goldman-Millson's theory, the kernel $K$ as defined in the theorem is exactly the image of this map. 
The map is a map of CMHS (the target is even pure) so the image is a sub-CMHS which is condition (3).

We just have to check condition (4). 
For $F$, the strictness is a consequence of the fact that the Goldman-Millson formality isomorphism trivializes the 
Hodge filtration. In other words, the $\C^{\ast}$ action on the algebra of forms $(\ker (D'),D'')$
gives a $\C^{\ast}$ action on the local ring $\Orn$ inducing the filtration.  
The map $\mu ^n$ preserves the decomposition, so it strictly preserves the filtration.

The statement for $G$ follows from the same argument.
Finally, for $W$ the statement can be seen by using the GM isomorphism $Spec(A)\cong (\Omega _{\rho}\times T)_n$.
\end{prf}

\begin{rem}
This MHS on $\Orn$ is uniquely defined and depends on only the base point $x\in X$. 
\end{rem}

\subsection{The universal VMHS}

\begin{theo}
 The universal representation $\rho_n:\pi_1(X,x)\to GL(\V_{\rho,x}\otimes \Orn)$ is the monodromy of  a 
$Gr$ polarizable $\C$-VMHS.
\end{theo}

\begin{rem}
This VMHS is uniquely defined and depends on only the base point $x\in X$. 
\end{rem}

\begin{prf}
On the bundle $\V_{\rho}\otimes \Orn$ over $X$, we have filtrations $sF$, $sG$ and $sW$ coming from the splitting of the Hodge decomposition 
on $V$, and the given filtrations on $\Orn$. We also have a universal Maurer Cartan element $(\eta , e^r)$ where $e^r$ can be viewed as a framing at the point $x\in X$. However, the universal Maurer Cartan element need not compatible
with the filtrations. Hence, there are gauge transformations $e^f$, $e^g$ and $e^w$ such that $e^f\cdot (\eta , e^r)$ is compatible with $sF$,
$e^g\cdot (\eta , e^r)$ is compatible with $sG$, and
$e^w\cdot (\eta , e^r)$ is compatible with $sW$. Define 
$$
F:= e^{-f}(sF), \;\;\; 
G:= e^{-g}(sG), \;\;\; 
W:= e^{-w}(sW). 
$$
These give three filtrations on the bundle $V\otimes \Orn$. The connection $\nabla + \eta$ is then compatible with these, in the sense
of Griffiths transversality for $F$, anti-transversality for $G$, and preserves $W$ (these are because of how the filtrations $F$, $G$ and $W$
were defined on the algebra of forms).  The action of $\Orn$ preserves these filtrations, and it follows from an analogue of the argument of
the previous section that these filtrations define MHS's at each point; so we get a VMHS. Furthermore, the weight-graded pieces are polarizable. 
\end{prf}

We conclude this section stating two obvious properties of the above constructions 
\begin{lem}
The natural maps $\mathcal O_{\rho | n+1} \to \Orn$ are morphisms of Mixed Hodge Artin local rings. 
\end{lem}

\begin{lem}
 We have a morphism of Mixed Hodge rings $\Orn \to End(\V_{\rho_n})$.
\end{lem}

\begin{rem}
The construction is independant of the K\"ahler form and functorial under morphisms $(Y,y)\to (X,x)$.
\end{rem}

 In order to see this, one has to adapt the argument of section \ref{vmhs} starting with the obvious $\alpha_1$
which is the tensor corresponding to the natural map $$H^1((X,x), End(\V_{\rho}) \to \mathcal{H}^1(X,End(\V_{\rho}))$$ given by composition of the Hodge isomorphism and the natural map $$H^1((X,x), End(\V_{\rho}) )\to H^1(X,End(\V_{\rho})).$$

\section{Proof of Theorem \ref{gengp}}

A more precise form of theorem \ref{gengp} is:

\begin{theo} \label{gengpb}
Let $G$ be a reductive algebraic group defined over $\C$. Let $\sigma: \Gamma\to G(\C)$ be a semisimple representation 
whose associated Higgs bundle is a fixed point of the $\C^*$-action on $M_{Dol}(X,G)$ \cite{Sim3}.

Let $\repg/\C$ be the affine scheme parametrizing the representations of $\Gamma$ with values in $G$ endowed with 
the action of $G$ by conjugation. Let $\hat\Omega_{\sigma}$ be the formal germ at $[\sigma]$ of the orbit of $\sigma$. 
There is a preferred isomorphism $GM^c: Spf(\Os) \to \hat\Omega_{\sigma}\times T$ where $T$ is the formal germ at the origin of the 
homogenous quadratic c\^one attached to the Goldman-Millson obstruction map $$S^2H^1(X, \mathrm{ad}_{\sigma})\to H^2(X,\mathrm{ad}_{\sigma}).$$

Let $\hat O_{\sigma}$ be the complete local algebra of $\repg$ at $\sigma$.
Then $\hat O_{\sigma}/\m^n$ carries a functorial $C$-MHS whose weight filtration comes from the powers of the ideal defining the orbit of $\sigma$. 

Let $\alpha$ be a rational  representation of $G$ with values in $GL_N$
 and let $\sigma_n:\Gamma\to G\Otn$ be the tautological representation defined in terms of $GM^c$. Denote by $\V_{\alpha,\sigma}$ the  local system in $\Otn$ free modules on $X$
 attached to the representation $\alpha\circ \sigma:\Gamma\to GL_N\Otn$.
  
  The $\C$-local system 
 underlying $\V_{\alpha,\sigma}$
 is the holonomy of a graded polarizable VMHS whose weight filtration is given by $$W_{-k} \V_{\alpha,\sigma}= \m^k.\V_{\alpha,\sigma} \quad k=0,..., n.$$ 
\end{theo}

In this section we give a proof of Theorem \ref{gengp} by spelling out  the differences between this more general case and the case $G=GL_n$ treated in Theorem \ref{gln}. 

Let $E$ be the real Zariski closure of the monodromy group of $\sigma$. This is a real reductive subgroup
of $G(\C)$ viewed as a real reductive group. $E$ is also of Hodge type \cite{Sim2}.

Recall that a real reductive algebraic group $E$ is said to be of
Hodge type if there is a morphism of real algebraic groups $h:U(1)\to Aut(E)$ such that $h(-1)$
is a Cartan involution of $E$, see \cite[p.46]{Sim2}.  By definition,
$h$ is a Hodge structure on $E$.  Connected groups of Hodge type are
precisely those admitting an isotropic Cartan subgroup. 

Consider a finite dimensional complex representation
 of $E$ 
$\alpha: E\to GL (\mathbb{V}_{\C})$. 
\begin{lem}$\ker(\alpha)$ is fixed by $h$. 
\end{lem}
\begin{prf}
By \cite[p.63, proof of lemma 5.5]{Sim2}
there is a morphism of real algebraic groups $z:U(1) \to E$ and an isogeny $\pi:U(1)\to U(1)$
such that $\mathrm{ad} (z)= h(\pi(z))$. 
Let $z',z\in U(1)$ such that $\pi(z)=z'$. Then $h(z')= \mathrm{ad}(z)$. Hence $h(z').g \in \ker(\alpha)$ iff 
 $zgz^{-1} \in \ker(\alpha)$ iff $g\in \ker(\alpha)$. 
\end{prf}

By \cite[lemma 5.5]{Sim2}, $\mathbb{V}_{\C}$ inherits a pure polarizable Hodge
structures of weight zero and $\alpha$ is a Hodge representation in the sense of \cite[lemma 5.6]{Sim2}. 
Hence, the local system attached to $\alpha\circ \sigma$ underlies a polarized $\C$-VHS. There is no uniqueness
since the polarization is 
not uniquely defined and the Hodge filtration can be shifted by an integer depending on each irreducible 
component of $\alpha$. 

The adjoint representation of $G(\C)$ restricted to $E$ gives rise to a $\C$-VHS on $X$ of weight zero which we call $\mathrm{ad}_{\sigma}$
If we keep track of the real structures then we can eliminate the shift of the Hodge 
filtration as a source of non-uniqueness but not the polarization.

Then, we can rewrite the construction of subsection \ref{gm} replacing $GL_N$ by $G$, using the new definition for $H^{\hdot}$ and $E^{\hdot}$
 given by $H^{\hdot}=H^{\hdot}(X, \mathrm{ad}_{\sigma})$ and $E^{\hdot}=E^{\hdot}(X, \mathrm{ad}_{\sigma})$
and the new augmentation vith values in $\mathfrak g$ defined by evaluation at $x$. 

To adapt subsection \ref{mhs}, we use the fact that $p=[-;-]: S^2H^1(X, \mathrm{ad}_{\sigma})\to H^2(X, \mathrm{ad}_{\sigma})$ respects the Hodge structure to  deduce that
$$I_2= \mathrm{Im}(p^t:H^2(X, \mathrm{ad}_{\sigma} )^*\to S^2 H^1(X,  \mathrm{ad}_{\sigma})^*)$$ is a weight $-2$ $\C$-HS. Hence, 
$$I_n= I_2 S^{n-2}H^1(X,  \mathrm{ad}_{\sigma})^*\subset S^n H^1(X,  \mathrm{ad}_{\sigma})^*$$ 
is a weight $-n$ sub $\C$-Hodge structure. 
Hence, 
$$\Pi_n = S^n H^1(X,  \mathrm{ad}_{\sigma})^*/I_n $$ 
is a weight $-n$ $\C$-HS and the algebra $S_2^H:=\sum_{n\in \N} \Pi_n$  viewed as an infinite dimensional split $\C$-MHS is naturally endowed with an algebra structure repsecting the MHS. The weight filtration
is exactly $W_{-n} S_2^H= \m_2 ^n$ since $\m_2= \sum_{n>0} \Pi_n$ is a maximal ideal of the complete local algebra $S_2^H$ and $\Pi_1$ generates $\m_2$ . Then we have $\Os\simeq S_2^H$. 

Section \ref{vmhs} is also easily adapted by replacing using a basis of $\mathcal{H}^1(X, \mathrm{ad}_{\sigma})$ in place of a basis 
of $\mathcal{H}^1(X,\mathrm{End}(\V_{\rho})$ in lemma \ref{a1} and the rest of the argument goes through without any difficulty.

In particular, we can put yet another MHS on the Artin local ring $\hat{O}_{\sigma|n}$ which corresponds to the $n$-th 
infinitesimal neighborhood of $\sigma$ and  whose weight filtration
is given by the powers of the ideal defining the orbit of $\sigma$. 
We can also interpret 
 the universal representation $\pi_1(X,x)\to G(\hat{O}_{\sigma|n})$ as the monodromy of  a 
$Gr$ polarizable $\C$-VMHS. 

This concludes the proof of theorem \ref{gengpb}.

There are several natural properties of the present construction
that we have not fully developped yet. For instance, the MHS on $\hat{O}_{\sigma|n}$  is likely to vary in a $\C$-VMHS when $x$ varies. 
We hope to understand this in a future work which should 
 do at the same time the comparison with \cite{Ha2}.

\end{document}